\documentclass[final,10pt,journal,letterpaper]{IEEEtran}
\pdfoutput=1
\IEEEoverridecommandlockouts
\usepackage{tabularx,ragged2e,booktabs}
\usepackage{stmaryrd}
\usepackage{arydshln}
\usepackage{graphicx}
\usepackage{amsmath}
\usepackage{amssymb}
\usepackage{mathrsfs}
\usepackage{url}
\usepackage[compress]{cite}
\usepackage{mdwlist}
\usepackage{paralist}
\usepackage{listings}
\usepackage{algorithm}
\usepackage{tikz}
\usepackage{pgfplots}
\usetikzlibrary{pgfplots.groupplots}
\allowdisplaybreaks 
\usepackage{tkz-graph}
\usepackage{multirow}
\usepackage{color}
\usepackage{bbm}
\usepackage[nomargin,inline,index]{fixme}
\usepackage{mathtools}
\usepackage{epstopdf}
\usepackage{enumitem}
\usepackage{dsfont}
\usepackage{tikzscale}
\usepackage{verbatim}
	\definecolor{mainblue}{rgb}{0.0, 0.0, 1.0}
\usetikzlibrary{arrows.meta}
\pgfplotsset{width=7cm,compat=1.8}
\newcommand{\al}[1]{{\color{black}#1}}

\usepackage{mathrsfs}
\usepackage{algpseudocode}
\usepackage{mathtools, cuted}
\usepackage{todonotes}
\usepackage{bm}

\newcolumntype{L}{>{\raggedright\arraybackslash}X}
\colorlet{pink}{red!40}
\colorlet{blue}{cyan!60}
\DeclarePairedDelimiter{\ceil}{\lceil}{\rceil}

\PassOptionsToPackage{usenames,dvipsnames,svgnames}{xcolor}
\usetikzlibrary{shapes.geometric,shapes.misc,arrows,positioning,automata,calc}

\newtheorem{lemma}{\sc Lemma}[section]
\newtheorem{theorem}[lemma]{\sc Theorem}

\newtheorem{proposition}[lemma]{\sc Proposition}

\newtheorem{remark}{\sc Remark}[section]
\newtheorem{assumption}{\sc Assumption}[section]

\newtheorem{problem}{\sc Problem}
\definecolor{qqqqff}{rgb}{0,0,1}
\definecolor{ffwwqq}{rgb}{1,0.4,0}
\definecolor{wqwqwq}{rgb}{0.3764705882352941,0.3764705882352941,0.3764705882352941}
\definecolor{dtsfsf}{rgb}{0.8274509803921568,0.1843137254901961,0.1843137254901961}
\definecolor{ffqqqq}{rgb}{1,0,0}
\definecolor{rvwvcq}{rgb}{0.08235294117647059,0.396078431372549,0.7529411764705882}
\definecolor{wvvxds}{rgb}{0.396078431372549,0.3411764705882353,0.8235294117647058}

\tikzstyle{block} = [draw, fill=black!20, rectangle,
    minimum height=3em, minimum width=3em]
\tikzstyle{sum} = [draw, fill=blue!20, circle, node distance=2cm]
\tikzstyle{input} = [coordinate]
\tikzstyle{output} = [coordinate]
\tikzstyle{pinstyle} = [pin edge={to-,thin,black}]
\tikzstyle{startstop} = [rectangle, rounded corners, minimum width=2cm, minimum height=0.5cm,text centered, draw=black, fill=red!30]
\tikzstyle{io} = [trapezium, trapezium left angle=70, trapezium right angle=110, minimum width=2cm, minimum height=0.5cm, text centered, draw=black, fill=blue!30]
\tikzstyle{process} = [rectangle, minimum width=3cm, minimum height=0.5cm,text width=3cm, text centered, draw=black, fill=orange!30]
\tikzstyle{decision} = [diamond, minimum width=1cm, minimum height=1cm, text centered, draw=black]%fill=green!30

\DeclareMathOperator*{\argmax}{arg\,max}

\renewcommand{\footnoterule}{%
  \kern -2pt
  \hrule width 0.3\textwidth height .5pt
  \kern 2pt
}

\begin{document}

% \author{Adair~Lang,~Farhad~Farokhi,~Michael~Cantoni, and~Iman~Shames\vspace{-.2in}
% \thanks{The authors are with the Department of Electrical and Electronic Engineering, University of Melbourne, Parkville, Victoria 3010, Australia. Emails:\{ffarokhi,cantoni,ishames\}@unimelb.edu.au} }

\title{
\al{Rigid-profile input scheduling under constrained dynamics with a water network application}%Load scheduling under constrained linear dynamics
  \thanks{
%   {Draft Manuscript: v05 November 10, 2017, v06 February 15,
%     2018, v07 August 17, 2018, v08 November 2, 2018 v09 November 25,
%     2018 v10 January 24, 2019 v11 February 6, 2019 v12 February 6,
%     2019, v13 August 5, 2019, v14 \today.
%This work extends on the authors work presented in \cite{Farokhi2016} and \cite{5400193}.
This work was supported by Rubicon Water~Pty.~Ltd., the Australian Research Council (LP130100605,~LP160100666), and a McKenzie Fellowship.}}

\author{Adair Lang, Michael Cantoni, \and Farhad Farokhi, Iman Shames
\thanks{A.~Lang is with Dept.~of Electrical and Electronic Engineering, The University of Melbourne, VIC 3010, Australia, and Rubicon Water, 1 Cato St, Hawthorn East, VIC 3123, Australia. {\tt\footnotesize aclang@student.unimelb.edu.au}}
\thanks{M.~Cantoni is with Dept.~of Electrical and Electronic Engineering, The University of Melbourne, VIC 3010, Australia. {\tt\footnotesize cantoni@unimelb.edu.au}}
\thanks{F.~Farokhi is with the Dept.~of Electrical and Electronic Engineering, The University of Melbourne, VIC 3010, Australia, {\tt\footnotesize ffarokhi@unimelb.edu.au}}
\thanks{I.~Shames is with Dept.~of Electrical and Electronic Engineering, The University of Melbourne, VIC 3010, Australia. {\tt\footnotesize ishames@unimelb.edu.au}}}

\maketitle

\begin{abstract}
  The motivation for this work stems from the problem of scheduling
  requests for flow at supply points along an automated network of
  open-water channels. The off-take flows are rigid-profile inputs to the system dynamics. In particular, the channel operator can only shift orders in time to satisfy constraints on the automatic response to changes in the load. This leads to a non-convex semi-infinite programming problem, with sum-separable cost that encodes the collective sensitivity of end users to scheduling delays. The constraints encode the linear time-invariant continuous-time dynamics and limits on the state across a \al{continuous} scheduling horizon. Discretization is used to arrive at a more manageable approximation of the semi-infinite program. A method 
  for parsimoniously refining the discretization is applied to ensure continuous-time feasibility for solutions of the approximate problem. It is then shown how to improve cost without loss of feasibility. Supporting analysis is provided, along with simulation results for a realistic irrigation channel setup to illustrate the approach.
  
\end{abstract}

\section{Introduction}
% Consider a continuous-time dynamical system that must satisfy hard
% constraints at all times over a finite time horizon. The system can
% admit a collection of pre-specified inputs termed \emph{loads}. These
% loads are rigid and can only be shifted back-and-forth in time to
% satisfy the constraints. The problem at hand is to determine a
% scheduling of these loads by optimizing an objective, such as
% minimizing scheduling delays, while ensuring constraint satisfaction
% at all times. This atypical scheduling problem can be formulated as a
% semi-infinite program (SIP), which is intractable in general. This
% paper explores computationally-friendly approximations of this problem
% that guarantee constraint satisfaction and provide a guaranteed level
% of sub-optimality.

The problem considered in this paper is motivated by
issues that arise in the management of automated open-channel networks for rural water
distribution~\cite{5400193}. Specifically, \al{the scheduling problem of interest pertains to} the timing of rigid-profile off-take flows at supply points, given constraints on the transient channel water-level and flow
responses to load changes, while minimizing the collective cost of
scheduling delays to end users. \al{Problems of this kind may also be relevant in other domains, e.g., energy systems and process control.}

Some of the most widely studied
scheduling problems, such as job and machine allocation,
typically involve only static relationships between the variables of
interest~\cite{pinedo2016scheduling, blazewicz1986scheduling,
  Allahverdi2015}. Similarly, within the literature on
irrigation networks, most scheduling studies are limited to static (or
steady-state) capacity constraint satisfaction, with no regard for the
transient behaviour~\cite{hong2012optimization, reddy1999optimal}.

\al{Scheduling problems that involve constraints on the transient behaviour of a dynamical system are considered in \cite{zachar2019,dias2018,allman2019,otashu2020,burnak2019}. 
In these works, the dynamics correspond to a discrete-time model, or a uniformly sampled continuous-time model, with constraint satisfaction only enforced on a uniformly sampled subset of the scheduling horizon. The issue of continuous-time constraint satisfaction is not addressed. Specific consideration of rigid-profile input scheduling subject to constraints on discrete-time dynamics is studied in \cite{5400193,Farokhi2016b}.}

\al{Irrigation channels are complex physical systems with continuous-time dynamics. For channels operating under closed-loop control, linear time-invariant models are suitable when there is ample in-channel storage \cite{WEYER2001,MichaelCSM}.
%In order to approximate the transient behaviour accurately  using discrete-time models with homogeneous sampling time  the sampling resolution would need to be fine enough. This would result in numerous constraints and huge models. 
This motivates the approach pursued below, where} attention is focused on optimizing scheduling decisions subject to constraints, over a continuous interval, on the evolution of a continuous-time model. 
%\al{across all of a continuous scheduling horizon}. 
In particular, a non-convex semi-infinite programming formulation of the scheduling problem is considered for linear time-invariant dynamics. %; such models
%represent the closed-loop behaviour of automated channels well in practice~\cite{mareels2005systems,MichaelCSM}. 
\al{Scheduling problems subject to constrained continuous-time dynamics are also considered in~\cite{flores-tlacuahuac2006,zhuge2012,Farokhi2016,tsay2019}. However, these works do not explicitly consider important aspects of computing feasible solutions for the semi-infinite program. This motivates the focus here on ensuring constraint satisfaction over the continuous scheduling horizon. Particular emphasis is given to the construction of direct discretizations with parsimonious non-uniform sampling of time and constraints such
that solutions of the corresponding approximate problem are continuous-time feasible. This is the aim of the proposed first stage of the approach. In the second and final stage}, the cost of the first-stage feasible schedule is improved by a sequential quadratic programming approximation of the original problem. \al{While it is possible to extend the approach to more general dynamics, the linear time-invariant context of this paper enables explicit characterization of the state and derivatives with respect to the scheduling variable, at any specified time for a given initial condition.}
% The approach considered below can potentially be extended to more general dynamical systems however, here linear time-invariant dynamics are considered since, in this case, the state and derivatives with respect to the scheduling variable, evaluated at specified time for a given initial condition can be explicitly determined.}

% Constraints on the evolution of a continuous-time model over a continuous interval are considered within the context of scheduling in \cite{}. Although these works  

% Off-take scheduling subject to constraints on the dynamic
% response is considered in \cite{5400193,Farokhi2016b,Farokhi2016}. Below, attention is focused on optimizing scheduling decisions
% subject to constraints on the evolution of a continuous-time model, as in~\cite{Farokhi2016}. In particular, a non-convex semi-infinite
% programming formulation of the problem is considered
% for linear time-invariant dynamics; such models
% represent the closed-loop behaviour of automated
% channels well in practice~\cite{mareels2005systems,MichaelCSM}. The difficult
% semi-infinite program is tackled in two stages. Direct discretizations
% of the decision space and the constraints are first constructed such
% that solutions of the corresponding approximate problem are also
% continuous-time feasible. In the second stage, the cost of this
% schedule is improved by a sequential quadratic programming
% approximation of the original problem. 

In the rest the introduction, the semi-infinite program is formulated. Aspects of the discretization-based approximations are then introduced. To conclude, the contributions are summarized and an outline is given for the paper.

\subsection{Problem formulation}
\label{sub:linear_time_invariant_systems}
Consider a linear time-invariant system with dynamics 
\begin{align} \label{eqn:LTI}
\dot{x}(t)=Ax(t)+Bu(t)+\sum_{j = 1}^m E_j v_j (t - \tau_j),~x(0)=x_0.
\end{align}
In this model, $x(t)\in\mathbb{R}^{n_x}$ is the state,
$u(t)\in\mathbb{R}^{n_u}$ is the control input, and each
$v_j(t-\tau_j)\in\mathbb{R}$ is a shifted load input at time
$t\in\mathbb{R}_{\geq 0}$. 
%\farhad{
%with $\tau_j$ denoting the scheduling delay}. 
Let $\mathcal{T}:=[0,T]$ be the scheduling
horizon and let $\mathbb{N}_{[a,b]}$ denote the integers between $a>0$ and $b\geq a$ inclusive. Given the admissible shift interval
$\mathcal{D}_j:= [\underline{\tau}_j, \overline{\tau}_j]$ for
each piecewise continuous load input request
$$v_j=(t\in[-\overline{\tau}_j,T-\underline{\tau}_j]\mapsto
v_j(t)\in \mathbb{R}), \quad~ j\in\mathbb{N}_{[1,m]},$$ 
the scheduling problem is to select the piecewise
continuous control input $u:\mathcal{T}\rightarrow\mathbb{R}^m$ and
load input schedule shift variables
%\farhad{delays} 
$\tau_j\in\mathcal{D}_j$
%,
%with $\overline{\tau}_j\geq \underline{\tau}_j$,
to satisfy linear constraints on the
%so that
%\begin{equation}
%  Cx(t;u,(v_j,\tau_j)_{j=1}^m)\leq c,~t\in [0,T],
%  \label{eqn:con_sol}
%\end{equation}
state 
\begin{subequations}
  \begin{align}
    &x(t;u,(v_j,\tau_j)_{j=1}^m)
    =
      \bar{x}_0^u(t)+\sum_{j=1}^m\bar{x}_j^v(t-\tau_j), \quad t\in\mathcal{T}, \label{eqn:xx}
    \intertext{where}
    &\bar{x}_0^u(t):=\exp(At)x_0+\int_0^t
                    \exp(A(t-s))Bu(s)\mathrm{d}s, \label{eqn:barx0} \\
    &\bar{x}_j^v(t):=\int_0^t
                    \exp(A(t-s))E_jv_j(s)\mathrm{d}s.
                    \label{eqn:barxv}
  \end{align}
\end{subequations}
% where $x(t)\in\mathbb{R}^{n_x}$ is the state,
%$u(t)\in\mathbb{R}^{n_u}$ is the control input,
%$v_j(t)\in\mathbb{R}$, $\forall j$, are \emph{a priori} known
%piece-wise continuous functions modelling loads, $\tau_j \in
%\mathcal{D}_j:=[\underline{\tau}_j, \overline{\tau}_j]$, $\forall j$,
%are scheduling delays with $\overline{\tau}_j \geq
%\underline{\tau}_j$ denoting the comfort bounds on scheduling delays,
%and $T$ is the scheduling horizon.
The best choice, among the feasible possibilities, minimizes
$f(\tau_1,\dots,\tau_m):=\sum_{j=1}^m h_j(\tau_j)$. The
differentiable $h_j:\mathcal{D}_j \to \mathbb{R}_{\geq 0}$ represents
the sensitivity of user $j$ to shifting its request \al{(e.g., production cost impacts of delayed delivery)}.
%The associated rigid-load scheduling problem is given by the following problem. 
\begin{problem}[Fixed-profile load scheduling]\label{prob:scheduling}%
% Given
% \begin{equation}\label{eq:ProblemData}
%   \Xi :=
%   \left(A,B,C,c,T,\left(h_j(\cdot),E_j,v_j(\cdot),
%         \underline{\tau}_j,\overline{\tau}_j\right)_{j=1}^m\right)
% \end{equation}
  Given set $\mathcal{U}$ of admissible control inputs, and 
  state constraint parameters $C\in\mathbb{R}^{n_c \times n_x}$
  and $c\in\mathbb{R}^{n_c}$, solve
\begin{subequations}\label{eqn:optim:1}
\begin{align}
f^* :=\min_{u,(\tau_j)_{j=1}^m} &\, f(\tau_1,\dots,\tau_m) \\
  \mathrm{s.t.}\hspace{.15in}&
                               %x(t; u, (\tau_j)_{j=1}^m)=\bar{x}_0^u(t)+\sum_{j=1}^m\bar{x}_j^v(t-\tau_j),\label{eqn:optim:1:sb}\\[-.45em]
                               %&
                                  Cx(t;u, (v_j,\tau_j)_{j=1}^m)\leq
                                  c,~ t\in\mathcal{T}, \label{eqn:optim:1:sc}\\
& \tau_j \in \mathcal{D}_j,~ j\in\mathbb{N}_{[1,m]}, \\
& u \in \mathcal{U}.
\end{align}
\end{subequations}
% and
% $\bar{x}_0^u(t)$ and $\bar{x}_j^v(t-\tau_j)$, $\forall j$, are,
% respectively, the solutions of the ordinary different equation
% \eqref{eqn:LTI} for the initial condition and for the $j$-th rigid
% load given by
% \begin{subequations}
% \begin{align}
% \bar{x}_0^u(t)&=\exp(At)x_0+\int_0^t \exp(A(t-\beta))Bu(\beta)\mathrm{d}\beta, \label{eqn:barx0} \\[-.4em]
% \bar{x}_j^v(t)&=\int_0^t \exp(A(t-\beta))E_jv_j(\beta)\mathrm{d}\beta, \forall j\in\mathbb{N}_{[1,m]}. \label{eqn:barxv}
% \end{align}
% \end{subequations}
%%
% such that $f(\tau_1,\dots,\tau_m):=\sum_{j=1}^m h_j(\tau_j)$ is minimised, where $h_j : \mathcal{D}_j \to \mathbb{R}$ are differentiable.
% Furthermore, subject to the requirement that $x(t)$ described by \eqref{eqn:con_sol} satisfies the following constraints
% \begin{align}
%   Cx(t)\leq c, \quad t\in\mathcal{T}:= [0,T], \label{eqn:Contr}
% \end{align}
% where $T$ is a real number denoting the planning horizon, $C\in \mathbb{R}^{n_c\times n_x}$, and $c\in \mathbb{R}^{n_c}$.
\end{problem}

Several factors make it challenging to solve \eqref{eqn:optim:1} in general:
\begin{enumerate}
	\item non-convexity of the constraint set with respect to the decision variables $(\tau_j)_{j=1}^m$;
	\item infinitely many constraint in \eqref{eqn:optim:1:sc}; and
	\item infinite dimensionality of the control input $u$.
\end{enumerate}
The first factor relates to the influence of $(\tau_j)_{j=1}^m$ on the
state $x(t;u,(v_j,\tau_j)_{j=1}^m)$ in the inequality constraints
\eqref{eqn:optim:1:sc}; see~\cite{Farokhi2016}. Non-convexity can make
it difficult to distinguish between potentially multiple local minima
and the global minimum.
%The second issue makes it difficult to even find a local
%minimum.
%for this problem since there are an infinite number of
%constraints to check; see~\eqref{eqn:optim:1:sc}.
With the second factor, the problem can thus be classified as a
non-convex semi-infinite program (SIP); see \cite{hettich1993semi} for
a comprehensive overview of such optimization problems. The third
factor relates to the complexity of solving constrained
continuous-time optimal control problems. Given the linear
dependence of (\ref{eqn:xx}) on $u$, this can be largely
overcome by restricting the set of admissible controls to be
finite-dimensional, \al{which is a commonly employed technique \cite{Betts2010}.}
%; e.g., $\mathcal{U}$ can be such that $u$ is
%a discrete-time signal passed through a zero-order hold
%(as done in \cite{Farokhi2016}), and then \eqref{eqn:optim:1} can be
%formulated as a SIP with infinitely-many constraints but finitely-many
%decision variables.
%Due to the linear dependence of (\ref{eqn:xx}) on $u$, such an
%approach
Doing so yields a problem in which the main difficulty rests with how
the state dynamics are affected by the decision variables $(\tau_j)_{j=1}^m$;
i.e., factors 1) and 2). As such, the rest of the paper is focused on
these difficulties, while simplifying the developments by assuming
that the control input is constant; i.e.,
$\mathcal{U}=\{u=(t\in\mathcal{T}\mapsto
u_0\in\mathbb{R}^{n_u})\}$. Indeed, $u$ is subsequently dropped as a decision
variable.
%and becomes a fixed part of the constraint
%\eqref{eqn:optim:1:sc} through $\bar{x}_0^u(t)$; see
%\eqref{eqn:barx0}.

% Most of the complexity from this issue can be overcome by restricting the admissible controls to be finite-dimensional e.g., if $\mathcal{U}$ is chosen such that $u$ is assumed to be a discrete-time signal passed through a zero-order hold, as done in \cite{Farokhi2016}, then \eqref{eqn:optim:1} can be formulated as a semi-infinite-dimensional optimization problem where there is an infinite number of constraints but finite dimensional decision variables. In this case, due to the linear dependence of $u(t)$ in \eqref{eqn:LTI}, this does not significantly add to the complexity of the scheduling problem, which lies with how the dynamics are affected by the shifts $\tau_j$. Hence the remainder of the paper will focus on the first two issues and to simplify developments from here on the control-input will be assumed to be constant i.e., $\mathcal{U}=\{u:\mathcal{T}\to \{u_0\}\}$ where $u_0 \in \mathbb{R}^{n_u}$ is predefined. In this case $u$ is dropped as a decision variable and just becomes a constant part of the constraint through $\bar{x}_0^u(t)$.

%The following is a standard assumption in the optimization literature.
%It ensures that Problem~\ref{prob:scheduling} is feasible with at
%least one schedule that is strictly feasible.
\begin{assumption}[Strict feasibility]\label{ass:feasibleOrig}
  There exists at least one schedule of load request shifts
  $(\tau_j)_{j=1}^m$, with $\tau_j \in \mathcal{D}_j$, such that
  $C x(t; u_0, (v_j,\tau_j)_{j=1}^m) - c < 0$ for all $t\in \mathcal{T}$.
  %; i.e., there exist a strictly feasible schedule.
\end{assumption}

Problem~\ref{prob:scheduling} could be reformulated in terms of a
switched hybrid model of the dynamics, in which the switching
instances correspond to the scheduling parameters $\tau_j$ in
\eqref{eqn:xx}. However, it appears that doing so leads to even more
challenging problems. Firstly, the set of switch states grows
exponentially in the number of users, since the evolution of the
channel dynamics depends on the combination in which requested load
profiles are applied. Secondly, the natural reformulation does not
match better studied forms of switched dynamics optimization
problems~\cite{Zu2004a,Zhu2015,Das2008a}. Specifically, the switching
sequence would not be pre-determined as in~\cite{Zu2004a}, this
reformulation does not decompose in the way considered
in~\cite{Zhu2015}, and it does not seem possible to comply with the
rigid profile requirements within the framework of~\cite{Das2008a},
where switching sequence and switching instants are decision
variables together.
% For example, a predefined switching sequence is often assumed and
% the aim is to determining the switching instants
% only~\cite{Zu2004a}. Also, decoupled problems are
% considered~\cite{Zhu2015}, in which an inner problem of determining
% the switching sequence, and an outer problem of determining the
% switching instants. Such a decoupling for
% Problem~\ref{prob:scheduling} does not seem
% possible. In~\cite{Das2008a}, the switching sequence and switching
% instants are variables to be set. However, the complexity of the
% appraoch grows with the number of switching states in the model,
% which as noted above explodes with the number $m$ of end
% users. Further, for the problem at hand the switching sequence must
% comply with constraints that maintain rigidity of the load
% profile. It is not apparent how existing methods could be adapted to
% tractably encode such constraints.
As such, switched system models are not considered
further.
% Additionally, it is not clear how explict constraint satisfication
% ~\ref{eqn:optim:1:sc} can be given for this p
Instead, direct discretization approximations are pursued for
obtaining a first-stage feasible solution, used to initialize a
continuous variable local approximation method for the original SIP.

\subsection{Direct discretization-based approximations}

A common approach to SIPs is to relax the problem by sampling the constraints across the time
horizon~\cite{hettich1993semi}. Even so, Problem~\ref{prob:scheduling}
remains difficult. Hence related work has also considered
discretization of the decision spaces to yield an integer linear
program~\cite{5400193,Farokhi2016b}. Whilst still difficult, standard
methods and solvers exist for such problems; e.g., see~\cite{wolsey1998integer}. On the other hand, shortcomings of this existing work relate to the use of uniform
discretization, which can lead to \al{unnecessarily} large problems \al{in terms of the number of decision variables and number of constraints}, and \al{the lack of} guarantees on continuous-time feasibility of the solutions
obtained. Aspects of the discretization approach are now presented as
a precursor to a summary of the main contributions made in this paper.
%Indeed, these include a method for
%parsimonously refining the discretization parameters until
%continuous-time constraint satisfaction is achieved by the solution to
%the approximate problem.

Replacing $\mathcal{D}_j$
with the finite discretized decision space
\begin{align}
  \hat{\mathcal{D}}_j &:=
                        \{\tau_j^{(1)},\dots,\tau_j^{(N_j)}\}\subset
                        \mathcal{D}_j
                        \label{eq:dDSpace}
\end{align}
for $j\in \mathbb{N}_{[1,m]}$ yields the problem
\begin{subequations}\label{eqn:optim:dd}
\begin{align}
  \hat{f}^* := \min_{(\tau_j)_{j=1}^m} &\, f(\tau_1,\ldots,\tau_m) \\
               %\sum_{j=1}^m h_j(\tau_j),\\
  \mathrm{s.t.}\hspace{.15in}& Cx(t;(\tau_j)_{j=1}^m)\leq c, ~t\in \mathcal{T}, \label{eqn:optim:dd:sc}\\
& \tau_j \in  \hat{\mathcal{D}}_j \text{ for } j\in\mathbb{N}_{[1,m]},
\end{align}
\end{subequations}
where explicit dependence of $x$ on the fixed $u_0$ and given requests $v_j$ has been dropped for convenience. 
% $N_j
%\in \mathbb{Z}$ is the cardinality of the set $\hat{\mathcal{D}}_j$
%i.e., $\left|\hat{\mathcal{D}}_j\right|=N_j$.
In principle, problem \eqref{eqn:optim:dd} could be solved by exhaustive search as the decision space is finite. A more sophisticated approach is presented in
Section~\ref{sec:stage1}. The approach is based on iteratively solving the following related finite-dimensional optimization problems. 

The optimization problem \eqref{eqn:optim:dd} is a restriction of \eqref{eqn:optim:1}, and thus,
$\hat{f}^* \geq f^*$. In fact, it is possible that \eqref{eqn:optim:dd} is
infeasible for instances of the finite decision sets
$(\hat{\mathcal{D}}_j)_{j=1}^m$. Under
Assumption~\ref{ass:feasibleOrig}, this can be overcome by adding 
elements to these sets. Parsimonious augmentation is desirable as the complexity of solving \eqref{eqn:optim:dd} grows with the cardinality $N_j$ of $\hat{\mathcal{D}}_j$, $j\in\mathbb{N}_{[1,m]}$. A corresponding method is proposed in Section~\ref{sub:EHSIPA}.

Replacing $\mathcal{T}$ in \eqref{eqn:optim:dd} with a finite subset of sample times at which the state constraints must hold yields a relaxation. With
\begin{equation}\label{eq:dTSpace}
	\hat{\mathcal{T}}_i :=
        \{t_i^{(1)},\dots,t_i^{(T_i)}\}\subset\mathcal{T},\quad i \in \mathbb{N}_{[1,n_c]},
\end{equation}
the problem becomes
\begin{subequations}\label{eq:LBD}
\begin{align}
  \hat{f}^{*,L} :=\min_{(\tau_j)_{j=1}^m} & 
              \,f(\tau_1,\ldots,\tau_m)  \label{eq:LBD:min}\\
  %\sum_{j=1}^m h_j(\tau_j),\nonumber\\
  \mathrm{s.t.} \quad & C_ix(t; (\tau_j)_{j=1}^m)\leq c_i,
                    ~ t\in \hat{\mathcal{T}}_i, 
                    ~ i \in \mathbb{N}_{[1,n_c]}, \label{eq:LBD:sc} \\
            & \tau_j \in  \hat{\mathcal{D}}_j, ~ j\in\mathbb{N}_{[1,m]},
\end{align}
\end{subequations}
where $C_i$ and $c_i$ denote the $i$-th row of $C$ and $i$-th entry of $c$, respectively. The cardinality of each constraint discretization set $\hat{\mathcal{T}}_i$ is denoted by $T_i$, which may be zero. Further, $\hat{f}^{*,L} \leq \hat{f}^*$, by definition. 
%On the other hand, $\hat{f}^{*,L}$ could be greater than or less than $f^*$.

Problem \eqref{eq:LBD} is tractable in the sense that it can be transformed into an integer linear program. This is shown in
Section~\ref{ssub:mixed_integer_linear_program}. However,
%due
%to the relaxation of the constraints
it is possible for solutions of \eqref{eq:LBD} to violate \eqref{eqn:optim:dd:sc}, and thus, be infeasible for problem \eqref{eqn:optim:1}.  
%Further,
%$\hat{f}^{*,L} \leq \hat{f}^*$, by definition. 
In view of Assumption~\ref{ass:feasibleOrig}, and continuity of the state with respect to $t\in\mathbb{R}_{\geq 0}$, this issue 
%of solutions of \eqref{eq:LBD} with respect to %Definition~\ref{def:SIPFeasibility} 
can be overcome by adding 
sample times to the
sets $\hat{\mathcal{T}}_i\subset \mathcal{T}$ and tightening the constraint \eqref{eq:LBD:sc}.
%SIP-feasible as per Definition~\ref{def:SIPFeasibility}. 
Specifically, the tightened problem takes the form
\begin{subequations} \label{eq:UBD}
\begin{align}
  \hat{f}^{*,U} :=\min_{(\tau_j)_{j=1}^m}&\, f(\tau_1,\ldots,\tau_m)\\ %\sum_{j=1}^m h_j(\tau_j),\\
  \mathrm{s.t.} \hspace{.15in} & C_ix(t; (\tau_j)_{j=1}^m)\leq c_i -\epsilon^{g}, ~ t\in \hat{\mathcal{T}}_i,\; i \in \mathbb{N}_{[1,n_c]},\\
  & \tau_j \in  \hat{\mathcal{D}}_j,~ j\in\mathbb{N}_{[1,m]},
\end{align}
\end{subequations}
with $\epsilon^g>0$. This is a restriction of \eqref{eq:LBD}, whereby $\hat{f}^{*,U} \geq \hat{f}^{*,L}$.

Problem (\ref{eq:UBD}) is neither a
restriction nor a relaxation of \eqref{eqn:optim:1}. However, as suggested above, suitable
constructions of $\hat{\mathcal{T}}_i$ and
$\epsilon^{g}$ can be made so that solutions of \eqref{eq:UBD} satisfy \eqref{eqn:optim:dd:sc}, and thus, \eqref{eqn:optim:1:sc}. In this case, $\hat{f}^{*,U}\geq \hat{f}^* \geq \hat{f}^{*,L}$. 
%Further, it is possible to ensure $\hat{f}^{*,U}-\hat{f}^*$ is %smaller than some specified tolerance. 
To this end, a recent algorithm for general SIPs from~\cite{Djelassi2017} is adapted to \eqref{eqn:optim:dd}, as detailed in Section~\ref{sec:stage1}.

\subsection{Contribution} % (fold)
   \label{ssub:contribution}
   The main contribution is a
   two-stage approach to Problem~\ref{prob:scheduling}; see \eqref{eqn:optim:1}. The first stage involves iterative construction of the discrete decision spaces $\hat{\mathcal{D}}_j$, constraint discretization sets $\hat{\mathcal{T}}_i$, and constraint restriction parameter $\epsilon^g$ in \eqref{eq:UBD}, from initial values. The resulting finite-dimensional
   problem yields a good schedule of shifts for which \eqref{eqn:optim:1:sc} holds; i.e., a schedule that is close in objective value to the optimal value in \eqref{eqn:optim:dd}, where only the decision space is discretized, and feasible for the original problem. The aim of the second stage is to improve this schedule while maintaining feasibility. To this end, the original continuous decision space is re-instated, and a sequential quadratic
   programming (SQP) method is used to approximately solve
   \eqref{eqn:optim:1}, initialized from the feasible schedule.

   The proposed two-stage approach goes beyond prior work in ensuring continuous-time feasibility via parsimoniously constructed discretizations. In particular, a method for general SIPs from~\cite{Djelassi2017} is adapted to the scheduling problem and extended to accommodate combined discretization of the decision space and the constraints. The resulting discretizations are not necessarily uniform, unlike the purely discrete-time formulations considered in~\cite{5400193,Farokhi2016b}. This can lead to smaller integer linear programming problems in the first-stage. In the second-stage, an SQP method from~\cite{Price1990,Price1992} is used to improve the cost without loss of feasibility. As in~\cite{Farokhi2016}, which develops a continuous-variable penalty-function based method that does not necessarily lead to improvement of the cost, analytic formulations of the derivatives required to implement the SQP method are possible by virtue of the linear time-invariance of the underlying dynamics.
%hettich1993semi,Reemtsen1998}
   The proposed algorithm is ultimately demonstrated on a non-trivial simulation example that is based on models of an Australian irrigation channel that are used operationally in the field and historically realistic demand profiles.

   % in the discrete-time, this approach addresses the issue of SIP-feasibility, question i and iii. This is done by applying a recent hybrid-SIP method from \cite{Djelassi2017} to the scheduling problem, bearing in mind the need to manage the size of the resulting integer program, question ii. Using the solution of this discrete-time approach as a feasible initial condition for a continuous-time approach provides a novel link between the two approaches. A sequential quadratic programming (SQP) method for SIPs is applied as a continuous-time algorithm, which has not been applied specifically to the scheduling problem before. It is shown that this approach yields guaranteed improvements on the initial feasible schedule whilst maintaining feasibility, question iv, which is in contrast to the continuous penalty methods in \cite{Farokhi2016}. The proposed algorithm is demonstrated on a non-trivial, realistic example from automated irrigation networks, using validated models of real irrigation pools, with an example requested demand profile that is consistent with that of these particular pools. 
   % The algorithm utilises an recent hybrid-SIP method from \cite{Djelassi2017} to address feasiblity with a focus of 
   % that is as good as or better than any possible schedule from . Thi
   %This algorithm utilises a recent hybrid-SIP method from \cite{Djelassi2017} to address the feasibility question (i). 
   % subsubsection contribution (end)   
   \subsection{Outline} % (fold)
   \label{ssub:outline}
   % subsubsection outline (end)
   The rest of the paper is organized as
   follows. Section~\ref{sec:stage1}
   includes discussion of the SIP discretization procedure from
   \cite{Djelassi2017} and the modifications made to ensure finite
   termination when it is applied to \eqref{eqn:optim:dd} from an initially infeasible discretization of the decision spaces. 
   %This is summarized in
   %Algorithm~\ref{alg:EHSIPA}. 
   The continuous-variable SQP approach to improving the first-stage feasible schedule is
   developed in Section~\ref{sec:stage2}. A
   formal characterization of the combined two-stage algorithm is
   presented in Section~\ref{sub:discussion}. Supporting analytical results are provided throughout. Numerical results based on
   non-trival irrigation channel scenarios are discussed in Section~\ref{sec:numerical_results}. The paper is concluded with
   discussion of future research directions in
   Section~\ref{sec:conclusion}.

\section{Stage 1 - Discretization with feasibility} % (fold)
\label{sec:stage1}
The Hybrid-SIP algorithm (HSIPA) from \cite{Djelassi2017}, which builds on~\cite{Mitsos2011,Mitsos2015}, provides a mechanism for generating constraint discretization sets $\hat{\mathcal{T}}_i$, and a constraint restriction parameter $\epsilon^g$, such that an optimizing schedule for \eqref{eq:UBD} is feasible for \eqref{eqn:optim:dd} with $\hat{f}^{*,U}-\hat{f}^*$ less than a specified tolerance $\epsilon^f>0$. 

\al{
The HSIPA requires the following input parameters:
\begin{itemize}
 \item initial constraint discretization sets $\hat{\mathcal{T}}_i\subset\mathcal{T}$, $i\in\mathbb{N}_{[1,n_c]}$;
	\item initial constraint restriction $\epsilon^g\in\mathbb{R}_{>0}$ in \eqref{eq:UBD}; and
	\item desired optimality tolerance $\epsilon^f\in\mathbb{R}_{>0}$ for problem \eqref{eqn:optim:dd}.
\end{itemize}
Given these inputs the HSIPA iteratively determines upper and lower bounds for $\hat{f}^*$. These bounds converge to within the specified tolerance $\epsilon^f$ of $\hat{f}^*$. The bounds are generated by successively solving iterations of the following three subproblems:
}	%The HSIPA involves successively solving iterations of the following three subproblems:
i) lower-bound problem \eqref{eq:LBD}; ii) upper-bound problem \eqref{eq:UBD}; and iii) refinement problem:
\begin{subequations}\label{eq:RES}
\begin{align}
	-\eta^* \!:=\!\! \min_{(\tau_j)_{j=1}^m,\eta} &\, -\eta	\\
\mathrm{s.t.}\hspace{.15in}& \sum_{j=1}^m h_j(\tau_j) - f^{R} \leq 0\label{eq:RESfREScon}\\
& C_ix(t; (\tau_j)_{j=1}^m)\leq c_i  - \eta,~ t\in \hat{\mathcal{T}}_i,~ i \in \mathbb{N}_{[1,n_c]}, \label{eq:RESCicon}\\
& \tau_j \in  \hat{\mathcal{D}}_j,~j\in\mathbb{N}_{[1,m]}, \text{ and } \eta\in\mathbb{R},
\end{align}
\end{subequations}
for \al{given} target objective value $f^R>0$. \al{The results are used to update the constraint discretization sets $(\hat{\mathcal{T}}_i)_{i=1}^{n_c}$, the restriction parameter $\epsilon^g$, and target objective $f^R$, as described subsequently. Iterations proceed}, until the upper and lower bounds obtained lie within $\epsilon^f$ of each other.
%The roles of \eqref{eq:LBD} and \eqref{eq:UBD} are described further in Sections~\ref{ssub:lower_bound} and \ref{ssub:upper_bound}, respectively. The optimization problem \eqref{eq:RES} relates to the largest constraint restriction that yields an objective value hat is no worse than the specified constant $f^{R}$. This problem plays role in refining lower and upper bound, as discussed further in Section~\ref{sub:refine_bounds}.
\al{Each solve of one of the three subproblems results in a candidate scheduling solution $(\tau_j)_{j=1}^m$. For given candidate solution $(\tau_j)_{j=1}^m$, the maximum level of constraint violation for each constraint $i \in \mathbb{N}_{[1,n_c]}$ is given by:
\begin{equation}\label{eq:LLP}
	g_i^*((\tau_j)_{j=1}^m)) := \max_{t\in\mathcal{T}} C_i x(t; (\tau_j)_{j=1}^m) - c_i.
\end{equation}
If $g_i^*((\tau_j)_{j=1}^m))>0$, then the candidate schedule is not feasible for \eqref{eqn:optim:dd}. In this case, a point in time that corresponds to this maximum level of constraint violation is added to the constraint discretization set. Specifically, given such a candidate scheduling solution $(\tau_j)_{j=1}^m$, the constraint discretization set is updated as follows:
\begin{equation}\label{eq:disUpdate}
	\hat{\mathcal{T}}_i \gets \begin{cases} \hat{\mathcal{T}}_i \cup \{t_i^*((\tau_j)_{j=1}^m)\} & \text{if } g_i^*((\tau_j)_{j=1}^m)) > 0 \\
	\hat{\mathcal{T}}_i & \text{otherwise}\end{cases},\vspace{-0.04in}
\end{equation}
where
\begin{equation}\label{eq:LLP_maximizers}
 	t_i^*((\tau_j)_{j=1}^m)\! \in\! \mathcal{T}_i^*((\tau_j)_{j=1}^m) := \argmax_{t\in\mathcal{T}} \left(C_i x(t; (\tau_j)_{j=1}^m)\! - \!c_i \right)
\end{equation}
is a corresponding maximizer of the maximum level of constraint violation.
}
%Iterative augmentation of the sets $\hat{\mathcal{T}}_i$, based on the results of successive solves, can lead to eventual feasibility of solutions to~\eqref{eq:UBD} with respect to~ problem~\eqref{eqn:optim:dd}. 
% As in~\cite[Algorithm 2]{Djelassi2017}, updates of the constraint discretization sets correspond to augmentation with times corresponding to maximum constraint violation across the original continuous scheduling horizon, using certain solutions of the lower-bound, upper-bound and refinement subproblems. %That is, after each subproblem is solved to yield a schedule $(\tau_j)_{j=1}^m$, the sets $\hat{\mathcal{T}}_i$ are updated as
% Specifically, given such a scheduling solution $(\tau_j)_{j=1}^m$, 
% \begin{equation}\label{eq:disUpdate}
% 	\hat{\mathcal{T}}_i \gets \begin{cases} \hat{\mathcal{T}}_i \cup \{t_i^*((\tau_j)_{j=1}^m)\} & \text{if } g_i^*((\tau_j)_{j=1}^m)) > 0 \\
% 	\hat{\mathcal{T}}_i & \text{otherwise}\end{cases},
% \end{equation}
% for each constraint $i \in \mathbb{N}_{[1,n_c]}$, where 
% \begin{equation}\label{eq:LLP}
% 	g_i^*((\tau_j)_{j=1}^m)) := \max_{t\in\mathcal{T}} C_i x(t; (\tau_j)_{j=1}^m) - c_i
% \end{equation}
% is the maximum level of constraint violation, and
% \begin{equation}\label{eq:LLP_maximizers}
%  	t_i^*((\tau_j)_{j=1}^m)\! \in\! \mathcal{T}_i^*((\tau_j)_{j=1}^m) := \argmax_{t\in\mathcal{T}} \left(C_i x(t; (\tau_j)_{j=1}^m)\! - \!c_i \right)
% \end{equation}
% is a corresponding maximizer. 
Feasibility of a schedule with respect to the infinite constraints \eqref{eqn:optim:1:sc} corresponds to
\begin{equation}\label{eqn:Feasble_Q}
	g_i^*((\tau_j)_{j=1}^m)) \leq 0,\quad i \in \mathbb{N}_{[1,n_c]}.
\end{equation}
\al{The update mechanisms of the constraint restriction parameter $\epsilon^{g}$ and target objective $f^R$ are discussed within the next section. The optimality tolerance $\epsilon^f$ should be chosen to reflect the desired optimality gap for solving \eqref{eqn:optim:dd}, which likely depends on the given setup. }

Algorithm 2 in \cite{Djelassi2017} is extended here in two ways. The first extension provides a way to deal with the finite sets $\hat{\mathcal{D}}_j$, and thus, a potentially infeasible problem~\eqref{eqn:optim:dd} for the given algorithm initialization.
% To this end, a parameter is added for use in ensuring finite-termination of the algorithm when \eqref{eqn:optim:dd} is initially infeasible.
The introduction of a new parameter enables interlacing of the HSIPA with an update procedure for the sets $\hat{\mathcal{D}}_j$ that achieves eventual feasibility of \eqref{eqn:optim:dd} from an initially infeasible discretization. The corresponding update procedure is presented in Section~\ref{sub:EHSIPA}. Secondly, the development as presented here elucidates the handling of a finitely indexed set of infinite constraints (i.e., $n_c>1$). 
%, with a remark given on how to extend to multiple constraints. 

\subsection{HSIPA for the scheduling problem}

Next, the notation in this paper is mapped to that used in \cite{Djelassi2017}, and the key algorithm parameters are identified. Aspects of the procedures associated with the aforementioned lower-bound, upper-bound and refinement problems are described in relation to generating the required constraint discretization sets $\hat{\mathcal{T}}_i$ and constraint restriction parameter $\epsilon^g$. Generation of the required discrete decision spaces $\hat{\mathcal{D}}_j$ is the topic of Section~\ref{sub:EHSIPA}. Modifications made here in adapting the HSIPA from~\cite{Djelassi2017} to the scheduling problem are highlighted below. 
 
\subsubsection{HSIPA notation and parameters}
%The notation in \cite{Djelassi2017} differs from that used here. 
Table~\ref{tab:HSIPA_notatation} relates the labels and notation used in this paper to those used in~\cite{Djelassi2017}. Notice, in particular, the multiplicity of constraint discretization sets here, one for each row of $C$ in \eqref{eqn:optim:dd:sc}, compared to one set in the formulation of~\cite{Djelassi2017}.
%A key difference pertains to the $i\in\mathbb{N}_{[1,n_c]} constraint sampling sets $\hat{\mathcal{T}}_i$, one for each scalar linear constraint on the system state, compared with the one sampled constraint index set $\mathcal{Y}^{D}$ in \cite{Djelassi2017}.
\begin{table}[]
    \centering
    \begin{tabularx}{\linewidth}{c|c|L}
    \cite{Djelassi2017}     &  This paper & Description \\\hline
     (SIP)    & \eqref{eqn:optim:dd} & SIP to be solved\\
     $f(\bm{x})$ & $f(\tau_1,\dots,\tau_m)$ & Objective function\\
     $f^*$    & $\hat{f}^*$ & Optimal objective value\\
     $\bm{x}$ & $(\tau_j)_{j=1}^m$ & Decision variable(s)\\
     $\mathcal{X}$ & $(\hat{\mathcal{D}}_j)_{j=1}^m$ & Decision space\\
     $g(\bm{x},\bm{y})$ & $Cx(t; (\tau_j)_{j=1}^m)- c$ & Inequality constraint(s)\\
     $\mathcal{Y}$ & $\mathcal{T}$ & Index set of infinite constraints\\
     $\mathcal{Y}^D$ & $\hat{\mathcal{T}}_i,~ i \in \mathbb{N}_{[1,n_c]}$ & Constraint discretization set(s)\\
     (LBD) & \eqref{eq:LBD} & Lower-bound sub-problem\\
     (UBD) & \eqref{eq:UBD} & Upper-bound sub-problem\\
     (RES) & \eqref{eq:RES} & Refinement sub-problem\\
     (LLP) & \eqref{eq:LLP} & Lower level program \\
     \hline\\
    \end{tabularx}
    \caption{A map of notation/labels in \cite{Djelassi2017} to the notation in this paper}
    \label{tab:HSIPA_notatation}
\end{table}
The HSIPA involves many parameters. \al{The parameters listed at start of Section \ref{sec:stage1} are important within the specific context of the subproblems \eqref{eq:LBD} and \eqref{eq:UBD}.
% , important parameters include the following:
% \begin{itemize}
% 	\item initial constraint discretization sets $\hat{\mathcal{T}}_i\subset\mathcal{T}$, $i\in\mathbb{N}_{[1,n_c]}$;
% 	\item initial constraint restriction $\epsilon^g\in\mathbb{R}_{>0}$ in \eqref{eq:UBD}; and
% 	\item desired optimality tolerance $\epsilon^f\in\mathbb{R}_{>0}$ for problem \eqref{eqn:optim:dd}.
% % 	\item $\epsilon^{LBP},\epsilon^{UBP},\epsilon^{RES}, \epsilon^{LLP}$, which are respectively optimality tolerances for lower-bound \eqref{eq:LBD}, upper-bound \eqref{eq:UBD}, refine-bound \eqref{eq:RES} and lower-level sub-problems \eqref{eq:LLP};
% % 	\item $r^g$ and $r^{LLP}$, which define how much to refine the constraint restriction and lower-level optimality tolerance within HSIPA, respectively.
% % 	\item $l_{max}$, the maximum number of consecutive updates of discretization of constraints allowed in refinement procedure before returning to lower-bound procedure.
% \end{itemize}
The following} parameter is introduced here to enable a decision space update procedure to be activated when \eqref{eqn:optim:dd} is infeasible:
\begin{itemize}
\item $k_r^{max}\in\mathbb{N}$, the maximum number of consecutive constraint restriction updates allowed in the upper-bound procedure.
\end{itemize}

\subsubsection{Lower-bound} % (fold)
\label{ssub:lower_bound}
 The lower-bound procedure of the HSIPA corresponds to Lines 2--11 of \cite[Algorithm 2]{Djelassi2017}. It pertains to solving the relaxed problem \eqref{eq:LBD} to determine a lower-bound $\hat{f}^{*,L}$ for $\hat{f}^*$. Given a feasible solution of \eqref{eq:LBD}, the constraint discretization sets $\hat{\mathcal{T}}_i$ are updated as per \eqref{eq:disUpdate}. This update, and those made in subsequent procedures of the HSIPA, ensure that successive runs of the lower-bound procedure result in a bound that converges to within a specified optimality tolerance of $\hat{f}^*$. To ensure finite termination of the HSIPA for an initially infeasible problem \eqref{eqn:optim:dd}, the lower-bound procedure is modified here to initially check if \eqref{eq:LBD} is feasible. When it is not, the procedure terminates without updating the sets $\hat{\mathcal{T}}_i$, after setting a flag that is used in the extension of Section~\ref{sub:EHSIPA} to trigger an update the discrete decision space $\hat{\mathcal{D}}_j$.
%  , the following check is added to \cite[Algorithm 2]{Djelassi2017} after line 2:
%  \begin{itemize}
%      \item \lstinline{if} \eqref{eq:LBD} \lstinline{is infeasible;} $\chi \gets $\lstinline{infeasible and terminate},
%  \end{itemize}
%  where $\chi$ is a flag used in the extension presented in Section~\ref{sub:EHSIPA} to trigger an update the discrete decision space $\hat{\mathcal{D}}_j$.

% subsubsection lower_bound (end)
\subsubsection{Upper-bound} % (fold)
\label{ssub:upper_bound}
Lines 12--24 of \cite[Algorithm 2]{Djelassi2017} constitute the upper-bound procedure of the HSIPA. This procedure relates to solving \eqref{eq:UBD}, which restricts the discretized constraints by $\epsilon^g$. If \eqref{eq:UBD} is infeasible, then the restriction is gradually reduced, in steps $\epsilon^g \gets \epsilon^g/r^g$ for specified $r_g>1$, until it becomes feasible. The modification made here, relative to~\cite{Djelassi2017}, is to limit the number of such steps to $k_r^{max}$, after which the upper-bound procedure terminates and the HSIPA returns to the lower-bound procedure, which could itself terminate as infeasible, triggering augmentation of the decision spaces $\hat{\mathcal{D}}_j$. When \eqref{eq:UBD} is feasible, and the resulting schedule satisfies \eqref{eqn:Feasble_Q}, the corresponding $\hat{f}^{*,U}$ is an upper-bound for $\hat{f}^*$; i.e., the infinite constraints \eqref{eqn:optim:dd:sc} are feasible. The upper-bound procedure then terminates, after a further step of $\epsilon^g$ reduction, and the HSIPA proceeds to the refinement procedure. Otherwise, when \eqref{eqn:Feasble_Q} does not hold, the constraint discretization sets $\hat{\mathcal{T}}_i$ are updated, as per \eqref{eq:disUpdate}, and the upper-bound procedure is repeated, including such updates, until an upper bound is found for $\hat{f}^*$. 
%Line 23 of \cite[Algorithm 2]{Djelassi2017} reduces the restriction when \eqref{eq:UBD} is infeasible and iterates back to the beginning of the upper-bound procedure. To prevent an infinite loop, in the case where \eqref{eqn:optim:dd} is infeasible, a counter $k$ is added, initially set to zero, and line 23 \cite[Algorithm 2]{Djelassi2017} is replaced with Algorithm~\ref{alg:UpperBound}. In particular, at most $k_r^{max}$ updates to the constraint restriction occur even if the \eqref{eqn:optim:dd} is infeasible. The modification returns to the lower-bound procedure, Line 2 \cite[Algorithm 2]{Djelassi2017}, if the maximum number of restrictions is exceeded. If the upper-bound procedure terminates in less than $k_r^{max}$ updates, then a valid upper-bound has been found, the counter $k$ can be reset to $0$ and the restriction $\epsilon^{g}$ can be reduced for subsequent runs. 
These updates, and those of successive subprocedures of the HSIPA, eventually lead to an upper-bound that lies within a specified tolerance of $\hat{f}^*$, as established in~\cite[Lemma 4]{Djelassi2017}. 
\subsubsection{Refinement procedure} % (fold)
\label{sub:refine_bounds}
The role of the refinement procedure is to improve both the upper- and lower-bounds whilst avoiding over-population of the sets $\hat{\mathcal{T}}_i$. The approach is adapted from \cite{Tsoukalas2011} in the HSIPA developed in~\cite{Djelassi2017}. The sub-procedure corresponds to Lines 26--42 in \cite[Algorithm 2]{Djelassi2017}. A target objective value, $f^R$, is selected as the mid-point between the current upper and lower bounds. The corresponding refinement problem \eqref{eq:RES} is then solved and updates are made on the basis of the outcome.
% where $\eta$ represents amount of restriction on constraints. Let $f^{RES}$ be chosen such that $f^{RES} > f(\tau_1,\dots,\tau_m)$ for any schedule $(\tau_j)_{j=1}^m$ such that $\tau_j \in  \hat{\mathcal{D}}_j, \,\forall j\in\mathbb{N}_{[1,m]}$. 
The following lemmas reveal how solving \eqref{eq:RES} may be used to improve the upper- or lower-bounds, $\hat{f}^{*,U}$ and $\hat{f}^{*,L}$, respectively, as detailed in the subsequent remarks. 
\begin{lemma}\label{lemm:RESLB}
If $\eta^*$ in \eqref{eq:RES} satisfies $\eta^*<0$, then $f^{R} < \hat{f}^*$. 
\end{lemma}
\begin{IEEEproof}
See Appendix~\ref{sec:proof_of_lemma_lemm:reslb}.
\end{IEEEproof}
\begin{lemma}\label{lemm:RESUB}
 If an optimal schedule $(\tau_j^*)_{j=1}^m$ for \eqref{eq:RES} satisfies \eqref{eqn:Feasble_Q}, then $f^{R} \geq \hat{f}^*$.
\end{lemma}
\begin{IEEEproof}
Since the resulting schedule satisfies \eqref{eqn:Feasble_Q}, it is feasible for \eqref{eqn:optim:dd}, and hence, $\hat{f}^* \leq f(\tau_1^*,\dots,\tau_m^*)$. By the constraint \eqref{eq:RESfREScon}, $f(\tau_1^*,\dots,\tau_m^*) \leq f^{R}$. As such, $\hat{f}^* \leq f^{R}$.
\end{IEEEproof}

% \begin{remark}
% Lemma~\ref{lemm:RESLB} gives conditions for generating an updated lower bound. If this occurs an increased $f^{RES}$ can be selected e.g., $f^{RES} = \frac{f^{UBD,*} + f^{LBD,*}}{2}$ and this Refinement procedure repeated. This yields a halving of the distance between the upper and lower-bounds for each consecutive run that satisfies conditions for Lemma~\ref{lemm:RESLB} and hence can greatly improve the overall computational time of HSIPA. 
% \end{remark}
\begin{remark}
The refinement procedure is repeated whilst the conditions of Lemma~\ref{lemm:RESLB} are satisfied, updating the lower-bound to $f^R$ accordingly for each such run. This can improve the overall computational time of the HSIPA as each run halves the difference between the upper and lower bound. 
%\farhad{[Not clear to me why? I mean the improvement in the overall computational time.]}
\end{remark}
\begin{remark}
Satisfaction of \eqref{eqn:Feasble_Q} in Lemma~\ref{lemm:RESUB} implies $\eta^* \in \mathbb{R}_{\geq 0}$;
%whereby 
%Hence, if the conditions for %Lemma~\ref{lemm:RESUB} are satisfied, 
%$\eta^*$ is 
which is thus an upper-bound for the smallest constraint restriction that retains feasibility. Updating $\epsilon^{g}$ accordingly yields improvement of $\hat{f}^{*,U}$ in subsequent runs of the upper-bound procedure, as described in \cite{Djelassi2017}.
 % in subsequent runs of the upper-bound procedure the restriction parameter $\epsilon^g$ should be smaller than $\eta^*$
% As outlined \cite[Algorithm 2]{Djelassi2017}in Therefore in order to improve the upper-bound in subsequent runs of the upper-bound procedure the restriction parameter $\epsilon^g$ should be smaller than $\eta^*$.
\end{remark}
\begin{remark}
It is possible that neither Lemma~\ref{lemm:RESLB} nor Lemma~\ref{lemm:RESUB} apply, leading to Line $36$ of \cite[Algorithm 2]{Djelassi2017}. In this case, the constraint discretization is updated, as per \eqref{eq:disUpdate}, and \eqref{eq:RES} is re-solved. To moderate growth in the cardinality of the constraint discretization sets $\hat{\mathcal{T}}_i$, this process terminates after a specified $l_{max}\in\mathbb{N}$ consecutive solves of \eqref{eq:RES}, and the HSIPA returns to the lower-bound procedure. 
\end{remark}

 In practice, the optimization problem \eqref{eq:RES} can only be solved to a specified tolerance, say $\epsilon^{RES}>0$. So as outlined in \cite[Proposition 1]{Djelassi2017} it is possible to get an inconclusive result. In this case the HSIPA returns to the lower-bound procedure. 

\begin{figure*}[htbp]
\includegraphics[width=0.9\linewidth]{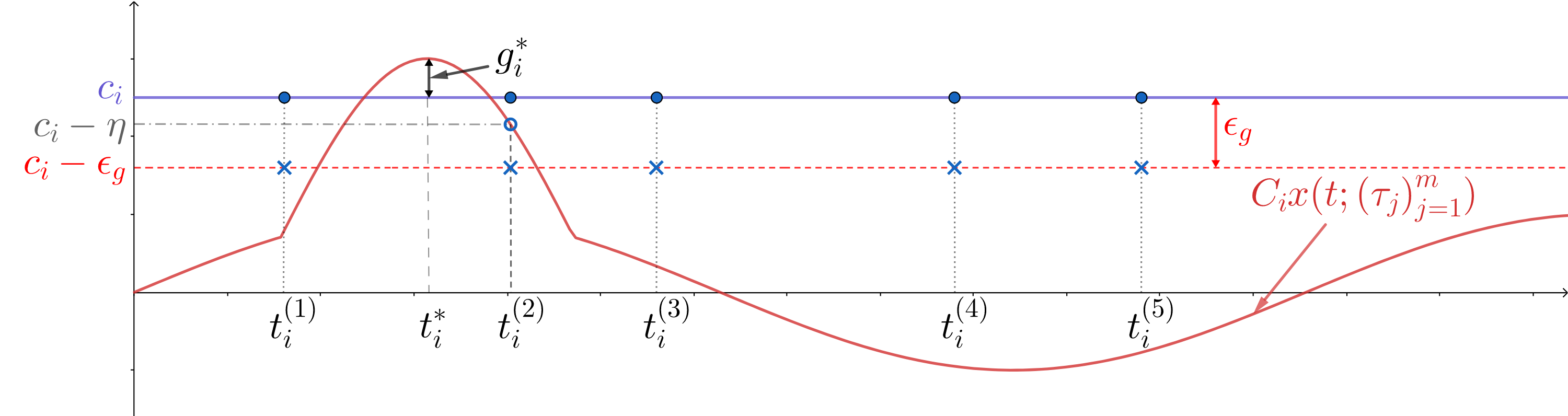}
\caption{\label{fig:hybridConstraintsEg}An example constraint trajectory to highlight the difference between the subproblems within the HSIPA.}
\end{figure*}
\subsubsection{HSIPA Illustration} % (fold)
\label{sec:hsipa_properties}
Fig.~\ref{fig:hybridConstraintsEg} shows an example trajectory of a constraint across a continuous scheduling horizon in order to illustrate the key differences between problems \eqref{eq:LBD}, \eqref{eq:UBD} and \eqref{eq:RES}. The original problem \eqref{eqn:optim:1} requires this trajectory to remain under the blue line ($c_i$) for the entire time-horizon. The lower bound problem~\eqref{eq:LBD} is feasible if the constraint trajectory lies below the blue dots at the specified sampled time-points, which holds true in this example. For the upper-bound problem \eqref{eq:UBD} the constraints are still enforced only at specified sample times, but the trajectory needs to be under the red-dashed line $c_i - \epsilon^g$. In this example, the constraint at $t_i^{(2)}$ is violated. When solving the refinement problem \eqref{eq:RES}, assuming this example trajectory has an objective value less than $f^{R}$, the maximum constraint restriction possible is $\eta$ shown as distance from $c_i$ to the open blue dot at $t_i^{(2)}$. In all cases, the problem is not SIP-feasible since $g_i^*$ is greater than 0. As such, the time sample $t_i^*$ would be added to the sampled time-horizon for this constraint for the next iterate.

\subsection{Extended HSIPA (EHSIPA)} % (fold)
\label{sub:EHSIPA}
Theorem 2 in \cite{Djelassi2017} establishes conditions for finite termination of the HSIPA at an $\epsilon^f$-optimal solution of the SIP at hand. For the result to apply here, the corresponding problem \eqref{eqn:optim:dd} must be feasible. This may not be the case for certain initializations of the decision space discretizations. The modifications made to the HSIPA presented in the previous section ensure finite termination, without a solution, when \eqref{eqn:optim:dd} is infeasible. In this case, the HSIPA sets a flag used to trigger augmentation of the decision space dicretizations. This continues until feasibility of \eqref{eqn:optim:dd} is achieved
. The subsequent run of the HISPA then augments the constraint discretization sets until the sub-problems \eqref{eq:UBD} yields an $\epsilon^f$-optimal solution of \eqref{eqn:optim:dd}, as described in the preceding subsection.
%In this section, the HSIPA is extended to include an update mechanism for the decision space discretization leading to eventual feasibility of \eqref{eqn:optim:dd}. 
%This is summarized in Algorithm~\ref{alg:EHSIPA}.

Given decision space discretizations $\hat{\mathcal{D}}_j = \{\tau_j^{(1)},\ldots,
\tau_j^{(N_j)}\}$, $N_j$ finite, $j\in\mathbb{N}_{[1,m]}$, the update mechanism
%in Algorithm~\ref{alg:EHSIPA} 
is defined by the following:
%At iteration $k$ with decision space $\hat{\mathcal{D}}_j[k]$, defined by \eqref{eq:dDSpace} with cardinality of the space $N_j[k]$, let the next decision space $\hat{\mathcal{D}}_j[k+1]$ be calculated as follows
\begin{align}
\hat{\mathcal{D}}_j \gets &
~\hat{\mathcal{D}}_j \cup
\left\{ \frac{\underline{\tau}_j + \tau_j^{(1)}}{2} 
\right\}
\cup 
\left\{ \frac{\tau_j^{(N_j)} + \overline{\tau}_j}{2}\right\} \nonumber \\
& \qquad\qquad \cup
\left\{ \frac{\tau_j^{(l-1)} + \tau_j^{(l)}}{2}~
:~l\in\mathbb{N}_{[2,N_j]}\right\} \label{eq:disDUpdate}.
\end{align}
% \begin{align}
% 	&\hat{\mathcal{D}}_j[k+1]  = \left\{\frac{\underline{\tau}_j +\tau_j^{(1)}}{2},\tau_j^{(1)},\frac{\tau_j^{(1)}+\tau_j^{(2)}}{2},\tau_j^{(2)},\dots,\right.\nonumber\\
% 	&\left. \tau_j^{(N_j[k]-1)},\frac{\tau_j^{(N_j[k]-1)}+\tau_j^{(N_j[k])}}{2},\tau_j^{(N_j[k])},\frac{\tau_j^{(N_{j}[k])}+\overline{\tau}_j}{2}\right\}\nonumber\\{}
% 	&\quad \forall j \in \mathbb{N}_{[1,m]}.\label{eq:disDUpdate}
% \end{align}
%The following lemma shows that this update mechanism eventually results in decision spaces that make \eqref{eqn:optim:dd} feasible.
% This leads to the following lemma as a way to eventually guarantee feasibility of the SIP on discretized decision spaces.
\begin{lemma}\label{lem:decisionHalve}
Under Assumption~\ref{ass:feasibleOrig}, a finite number of decision space discretization updates \eqref{eq:disDUpdate} 
% is such that 
%the decision spaces $\hat{\mathcal{D}}_j[k]$ for all $j\in\mathbb{N}_{[1,m]}$, i.e., $N_j[k]<\infty$, then 
% there exists a finite $\bar{k}\in\mathbb{N}$ such that 
leads to a feasibile problem
\eqref{eqn:optim:dd}. 
%is feasible for $\hat{\mathcal{D}}_j[\overline{k}]$, $j \in \mathbb{N}_{[1,m]}$.
% Assume Assumption~\ref{ass:feasibleOrig} holds, an initial decision space, $\hat{\mathcal{D}}_j^{(0)}$ is chosen with $N_j^{(0)}<\infty$ as defined in \eqref{eq:dDSpace}, for $j \in \mathbb{N}_{[1,m]}$ and assume that the discretization is iteratively updated as per \eqref{eq:disDUpdate} if \eqref{eqn:optim:dd} is found to be infeasible, then after a finite number, $\overline{k}$, of iterations the SIP \eqref{eqn:optim:dd} on discretized decision spaces $\hat{\mathcal{D}}_j^{(\overline{k})}$ will be feasible.
\end{lemma}
\begin{IEEEproof}
See Appendix~\ref{sec:proof_of_lemma_lem:decisionhalve}.
\end{IEEEproof}
% subsubsection feasibility (end)
%Lemma~\ref{lem:decisionHalve} guarantees eventual feasibility of the restricted problem \eqref{eqn:optim:dd}, requried for the main result presented in the subsequent sub-section.  Algorithm~\ref{alg:EHSIPA} shows the Extended HSIPA outlining how the update \eqref{eq:disDUpdate} is interlaced with the modified HSIPA.
\begin{remark}\label{rem:GrowthOfDecisionSpaces}
Each update \eqref{eq:disDUpdate} leads to a doubling of the cardinality of every $\hat{\mathcal{D}}_j$, and consequently larger integer linear programs to solve in the HSIPA. It is of interest to devise a more parsimonious mechanism, which is the topic of future work. Dependence of run time on the discretised decision space cardinalities is explored empirically in Section~\ref{sec:numerical_results}.

% The aforement , which requires no \emph{a prori} information, can potentially lead to extremely large integer programs to be solved as part of the HSIPA. This is because it results in a doubling of the number of samples $N_j^{(k)}$ at each iterate, i.e., $N_j^{(k+1)} = 2N_j^{(k)}$. If some information about the problem structure was known or could be easily determined it may be possible to reduce the number of samples added whilst guaranteeing eventual feasibility. This is a subject of future work.T
% Using \eqref{eq:disDUpdate}, which requires no \emph{a prori} information, may lead to many unnecessary samples being added to each $\hat{\mathcal{D}}_j^{(k+1)}$. Future work is to explore how to gain information about the problem or use \emph{a prori} information to reduce the number of samples added whilst guaranteeing eventual feasibility. The link between $N_j^{(k)}$ and run time is explored further in Sections~\ref{ssub:discretization_and_computational_time} and \ref{sec:numerical_results}.
\end{remark}

\begin{assumption}[Solution tolerances]\label{assum:NLP_Solver}
The feasibility of problems \eqref{eq:LBD} and \eqref{eq:UBD}
%\farhad{, the} %and \eqref{eq:RES}, 
is assessed without tolerance to infeasibility. When found to be feasible, these problems are solved to within specified global optimality tolerances 
%If any of the problems \eqref{eq:LBD}, \eqref{eq:UBD} and \eqref{eq:RES} are infeasible they are found to be infeasible when solving without any infeasibility tolerance. Any of these problems that are not infeasible and problem \eqref{eq:LLP} can be solved globally to within arbitrary optimality tolerances 
$\epsilon^{LBP},\epsilon^{UBP}
%\epsilon^{RES},
\in \mathbb{R}_{>0}$, respectively.
Further, it is possible to solve problem~\eqref{eq:LLP} to arbitrary optimality tolerance. 
\end{assumption}
\begin{remark}
In principle, \eqref{eq:LBD} and \eqref{eq:UBD} can be solved exactly by exhaustive search. However, this would be impractical for any reasonably sized problem. As shown in the next subsection, the discretization of decision spaces enables transcription to standard linear integer programs. While still difficult to solve, existing algorithms for such problems offer much greater efficiency than exhaustive search; see~\cite{wolsey1998integer,belotti2013}. Problem \eqref{eq:LLP} can be solved in principle by simulation of the linear time-invariant dynamics with sufficient accuracy.
\end{remark}

With Lemma~\ref{lem:decisionHalve}, the conditions required to apply~\cite[Theorem 2]{Djelassi2017} are met.
Its application yields the following result.
%reproduced below for completeness.
\begin{proposition}\label{th:HISPA}
%[Finite Convergence]
Given desired optimality tolerance $\epsilon^f$ for problem \eqref{eqn:optim:dd}, under Assumption~\ref{assum:NLP_Solver}, if there exists $\tilde{\epsilon}^f\in\mathbb{R}_{>0}$ such that 
$\tilde{\epsilon}^{f} \geq f(\tau_1^f,\dots,\tau_m^f) - \hat{f}^*$ for a schedule $(\tau_j^f)_{j=1}^m$ with 
$\max_{i\in\mathbb{N}_{[1,n_c]}} g_i^*((\tau_j^f)_{j=1}^m)) < 0$ (i.e., strictly feasible)
and
\begin{equation}\label{eq:tol_relations}
	\epsilon^{f} > \tilde{\epsilon}^{f} + \epsilon^{LBP} + \epsilon^{UBP},
\end{equation}
then the EHSIPA terminates after finite steps at a schedule $(\tau_j^*)_{j=1}^m$ that satisfies \eqref{eqn:optim:1:sc} and $f(\tau_1^*,\dots,\tau_m^*) \leq \hat{f}^*+\epsilon^{f}$.
% Let $\epsilon^{f}, \epsilon^{LBP},\epsilon^{UBP}\in\mathbb{R}_{>0}$ be set such that
% \begin{equation}\label{eq:tol_relations}
% 	\epsilon^{f} > \tilde{\epsilon}^{f} + \epsilon^{LBP} + \epsilon^{UBP},
% \end{equation}
% where $\tilde{\epsilon}^{f}$ satisfies $\tilde{\epsilon}^{f} \geq f(\tau_1^f,\dots,\tau_m^f) - \hat{f}^*$ for a schedule $(\tau_j^f)_{j=1}^m$ satisfying Assumption~\ref{ass:feasibleOrig}, Then, under Assumption~\ref{assum:NLP_Solver}, Algorithm~\ref{alg:EHSIPA} terminates finitely with a schedule $(\tau_j^*)_{j=1}^m$ that satisfies \eqref{eqn:optim:1:sc} and has an objective value $f(\tau_1^*,\dots,\tau_m^*)$ that satisfies $f(\tau_1^*,\dots,\tau_m^*) \leq \hat{f}^*+\tilde{\epsilon}^{f}$ i.e., an $\tilde{\epsilon}^{f}-$optimal schedule.
\end{proposition}
\begin{IEEEproof}
First note that $\mathcal{T}$ and $\hat{\mathcal{D}}_j$ are compact sets; the former is a bounded real interval and the latter a finite set~\cite{rudin1991functional}.
%, since every $\mathcal{S}$-cover of $\hat{\mathcal{D}}_j[k]$ has a finite sub-cover
Therefore, \cite[Assumption 1]{Djelassi2017} is satisfied. Now from \eqref{eqn:xx}, note that
the constraints are all continuous in the schedule variables, as is the objective function, since the $h_j$ are differentiable. Hence for the case of a single constraint, $n_c=1$, \cite[Assumption 2]{Djelassi2017} is satisfied. In the case of $n_c>1$, the multiple constraints can be replaced with a single constraint given by:
\begin{equation}\label{eq:maxCon}
	\max_{i\in\mathbb{N}_{[1,n_c]}} (C_i x(t; (\tau_j)_{j=1}^m) - c_i) \leq 0,~ t \in [0,T],
\end{equation} 
which is continuous since the maximum of a set of continuous functions is continuous. 
% The modifications made to the HSIPA in Section~\ref{sub:sample_time_and_constraint_} ensure that if \eqref{eqn:optim:dd} is infeasible the algorithm terminates in finite steps, the HSIPA, \cite[Algorithm 2]{Djelassi2017} will finitely terminate and therefore an update of discretization on Line 6 Algorithm~\ref{alg:EHSIPA} occurs. Therefore, Lemma~\ref{lem:decisionHalve} ensures \cite[Assumptions 3]{Djelassi2017} holds for some $\overline{k}$ and Assumption~\ref{assum:NLP_Solver} meets \cite[Assumption 4]{Djelassi2017}. 
Hence, under Assumption~\ref{assum:NLP_Solver}, all requirements for the proof in \cite{Djelassi2017} to follow are satisfied. 
\end{IEEEproof}

\begin{remark}
\al{Proposition~\ref{th:HISPA} only guarantees the number of steps required to reach optimality is finite. It does not give any guarantees on the size of this number. The chosen tolerance $\epsilon^f$ can have a significant effect on the overall number of steps required. Instead of letting the algorithm run until desired optimality guarantees are achieved, the algorithm can be stopped after at least one valid upper bound has been found. The schedule associated with the lowest upper-bound can then be returned, since this schedule is feasible, i.e., satisfies \eqref{eqn:optim:1:sc}. The sub-optimality gap achieved with such an approach would be highly problem and configuration dependent.}
\end{remark}

\subsection{Sub-problems as Integer Linear Programs} % (fold)
 \label{ssub:mixed_integer_linear_program}
\al{For the linear dynamics considered here,} the HSIPA sub-problems \eqref{eq:LBD}, \eqref{eq:UBD}, and \eqref{eq:RES}, can be equivalently reformulated as binary linear programs. To this end, define the following:
\begin{align*}
\Psi_{ji}&= 
 \begin{bmatrix}
 \bar{x}_j^{v(1,1)}          & \bar{x}_j^{v(1,2)} &\cdots &\bar{x}_j^{v(1,N_j)}  \\
 \bar{x}_j^{v(2,1)}          & \bar{x}_j^{v(2,2)} &\cdots  &  \bar{x}_j^{v(2,N_j)} \\
 \vdots                      & \vdots            &\ddots   &  \vdots  \\
 \bar{x}_j^{v(T_i,1)} & \bar{x}_j^{v(T_i,2)}      &\cdots &\bar{x}_j^{v(T_i,N_j)}
 \end{bmatrix}, 
\end{align*}
\begin{align*}
\bar{\Psi}_j&= \begin{bmatrix} (I_{T_1}\otimes C_1)\Psi_{j1} \\ (I_{T_2}\otimes C_2)\Psi_{j2} \\ \vdots \\ (I_{T_{n_c}}\otimes C_{n_c})\Psi_{jn_c} \end{bmatrix}, \qquad b_i=
\begin{bmatrix}
\bar{x}_0^u(t_i^{(1)}) \\ \bar{x}_0^u(t_i^{(2)}) \\ \vdots \\ \bar{x}_0^u(t_i^{(T_i))})\end{bmatrix},\\
	 \bar{c} &= \begin{bmatrix} (\mathds{1}_{T_1}\otimes {c}_1)-(I_{T_1} \otimes C_1)b_1 \\ (\mathds{1}_{T_2}\otimes {c}_2)-(I_{T_2} \otimes C_2)b_2\\ \vdots \\ (\mathds{1}_{T_{n_c}}\otimes {c}_{n_c})-(I_{T_{n_c}} \otimes C_{n_c})b_{n_c}\end{bmatrix}, \\
	 f_j&=\begin{bmatrix}
h_j(\tau_j^{(1)}) & h_j(\tau_j^{(2)}) & \cdots & h_j(\tau_j^{(N_j)})
\end{bmatrix}^\top,
\end{align*}
where $\bar{x}_j^{v(k,l)} = \bar{x}_j^{v}(t_i^{(k)}-\tau_j^{(l)})$ with $\tau_j^{(l)}$ and $t_i^{(k)}$ as defined in \eqref{eq:dDSpace} and \eqref{eq:dTSpace}, respectively.

Now the lower bound problem \eqref{eq:LBD} can  be rewritten as an equivalent binary linear program given by
\begin{subequations} \label{eqn:optim_central_3}
\begin{align}
\min_{(z_j)_{j=1}^m} & \sum_{j=1}^m f_j^\top z_j,\label{eqn:optim_central_3a}\\
\mathrm{s.t.}\hspace{.11in}
& \sum_{j=1}^m \bar{\Psi}_j z_j -\bar{c}\leq 0,\label{eqn:optim_central_3b}\\
&z_j\in \{0,1\}^{N_j},\quad \mathds{1}^\top z_j=1, \quad j\in\mathbb{N}_{[1,m]}.\label{eqn:optim_central_3c}
% & \mathds{1}^\top z_j=1, \forall j\in\mathbb{N}_{[1,m]}. \label{eqn:optim_central_3d}
\end{align}
\end{subequations}
Problem \eqref{eq:UBD} can be reformulated similarly, except $\bar{c}$ is formed with $c_i - \epsilon^g$ in place of $c_i$. Similarly, reformulation of \eqref{eq:RES} results in
\begin{subequations}\label{eqn:RES_D}
\begin{align}
 -\eta^* &= \min_{(z_j)_{j=1}^m, \eta} -\eta \\
\text{s.t } \quad & \sum_{j=1}^m f_j^\top z_j - f^{R} \leq 0\\
& \sum_{j=1}^m \bar{\Psi}_j z_j -\bar{c}\leq -\bar{\eta}, \\
&z_j\in \{0,1\}^{N_j}, \quad \mathds{1}^\top z_j=1, \quad j\in\mathbb{N}_{[1,m]},%\\
% & \mathds{1}^\top z_j=1, \forall j\in\mathbb{N}_{[1,m]},\nonumber
\end{align}
\end{subequations}
where $\bar{\eta} = \begin{bmatrix} (\mathds{1}_{T_1}\otimes {\eta})^\top & (\mathds{1}_{T_2}\otimes {\eta})^\top & \dots & (\mathds{1}_{T_{n_c}}\otimes {\eta})^\top \end{bmatrix}^\top$. 
% \eqref{eqn:RES_D} is a \emph{mixed}-integer linear program, due to the continuous variable $\eta$, however it can be solved by using similar algorithms to the completely integer linear programs.

The reformulations above generalize the approach in \cite{5400193,Farokhi2016b}, where purely discrete-time dynamics were considered. The generalization allows for the non-uniform sampling and different sampling sets for each constraint generated by the EHSIPA. Each reformulation involves a change of variables from $\tau_j$ to $z_j$ for each user $j$. The corresponding relationship between $z_j$ and $\tau_j$ is given by
\begin{equation}\label{eq:z_jCoV}
	\begin{bmatrix} \tau_j^{(1)} & \dots & \tau_j^{(N_j)} \end{bmatrix} z_j = \tau_j.
\end{equation}
This change of variable depends on $N_j, T_i < \infty$, $j\in \mathbb{N}_{[1,m]}$, $i \in \mathbb{N}_{[1,n_c]}$, which is not the case for the original formulation~\eqref{eqn:optim:1} where the corresponding sets are continuous intervals.

The costs and constraints in \eqref{eqn:optim_central_3} and \eqref{eqn:RES_D} are linear in the decision variables $(z_j)_{j=1}^m$, except for the binary requirement on the entries of $z_j$. Dedicated linear-integer programming solvers are readily available for such problems; e.g., see~\cite{wolsey1998integer}. The difficulty of solving discrete problems in non-standard forms, such as \eqref{eq:LBD}, \eqref{eq:UBD} and \eqref{eq:RES}, is highlighted in \cite{belotti2013}.

The reformulations \eqref{eqn:optim_central_3} and \eqref{eqn:RES_D} can be relaxed to linear programs by replacing the binary constraint with $z_j^{(q)} \in (0,1)$, $q\in \mathbb{N}_{[1,N_j]}$ and $j\in \mathbb{N}_{[1,m]}$. This relaxation leads to much nicer problems to solve, but at the cost of losing the rigidity constraint, as discussed in \cite{5400193}. So it is only applicable where the load-profile is not constrained to a fixed shape. 

\section{Stage 2 - Localized cost improvement}
\label{sec:stage2}
The EHSIPA described above involves discretization to arrive at a schedule $(\tau_j^*)_{j=1}^m$ that is feasible for \eqref{eqn:optim:1}. This schedule is $\epsilon^f$-optimal with respect to the restriction \eqref{eqn:optim:dd} associated with the final discretized decision spaces $\hat{\mathcal{D}}_j$. Since  $\hat{\mathcal{D}}_j\subset \mathcal{D}_j$ by definition, it may be possible to improve the sum of user sensitivities to scheduling delay, $f(\tau_1^*,\dots,\tau_m^*)$, by returning to the original continuous interval decision spaces. To this end, an SQP algorithm is used in a second stage of the proposed approach to Problem~\ref{prob:scheduling}. Starting with the feasible schedule generated by the EHSIPA, the approach leads to localised improvement of the cost whilst maintaining feasibility with respect to the original constraints in \eqref{eqn:optim:1}.

\subsection{SQP approximation} % (fold)
\label{sub:sequential_quadratic_program}
SQP methods involve solving quadratic programs formulated to approximate the original problem around each iterate~\cite{Nocedal2006}. Adaptations of such methods to SIPs, like  problem~\eqref{eqn:optim:1}, can be found in~\cite{hettich1993semi,Reemtsen1998}. The usual approach is to discretize the infinite constraints. However, such relaxation may result in solutions that are infeasible for the original problem. This could be overcome by refining the constraint discretization, in a similar fashion to the decision space update \eqref{eq:disDUpdate} in the EHSIPA, for example.
%, similar to how the decision space is refined in %\eqref{eq:disDUpdate}, 
%until the solution is SIP feasible. 
However, this can lead to a 
%large discretization and 
large number of constraints.
%, making the quadratic problem less tractable. 
%Alternatively, the constraints can be sampled at the 
The approach here is to sample constraint $i \in \mathbb{N}_{[1,n_c]}$ at the maximizers $\mathcal{T}_i^*((\tau_j)_{j=1}^{m})$ defined by \eqref{eq:LLP_maximizers} for the given feasible schedule iterate $(\tau_j)_{j=1}^m$. When the number of maximizers is finite, as subsequently assumed, a tractable SQP algorithm ensues~\cite{Price1990}. 
% The approach is summarized in Algorithm~\ref{alg:SQP}. The key components are outlined subsequently. 
% \begin{assumption}[Constraint lhs maximizers]\label{ass:FiniteMaximizers}
% For each constraint $i\in\mathbb{N}_{[1,n_c]}$, the cardinality $T_i^*:=|\mathcal{T}_i^*((\tau_j^*)_{j=1}^m)|$ of the maximisers in \eqref{eq:LLP} is finite; i.e., $\mathcal{T}_i^*((\tau_j^*)_{j=1}^m) = \{t_i^{*(1)},\ldots, t_i^{*(T_i^*)}\}$ with $T_i^*<\infty$.
% \end{assumption}

To construct the approximate quadratic program at each step, differentiability of the constraints is required. For the problem \eqref{eqn:optim:1} at hand, these derivatives can be formulated explicitly. 
\begin{lemma}\label{lem:constraintDef}
The left-hand side of constraint $C_i x(t;(\tau_j)_{j=1}^m)-c_i \leq 0$, $i\in \mathbb{N}_{[1,n_c]}$, in \eqref{eqn:optim:1} is differentiable in the schedule $(\tau_j)_{j=1}^m$. Moreover,
\begin{align*}
\frac{\partial}{\partial \tau_\ell}  C_i x(t;(\tau_j)_{j=1}^m)
%&=
%- C_i E_\ell v_\ell(t - \tau_\ell)  
%- C_i \int_0^{t- \tau_\ell} \!\!\!\!\!\!\! A \exp (A (t - \tau_\ell - \beta)) E_\ell v_\ell(\beta) d \beta\\
&= - C_i E_\ell v_\ell(t - \tau_\ell) 
 - C_i A \al{\bar{x}_\ell^v}(t-\tau_\ell)
\end{align*}
\begin{IEEEproof}
See Appendix~\ref{proof:lem:contraintDef}.
\end{IEEEproof}
\end{lemma}
For given feasible schedule $(\tau_j)_{j=1}^m$, the quadratic program to solve at each step of the SQP algorithm takes the form
\begin{subequations} \label{eqn:local_Quad_approx}
\begin{align}
\tilde{f}^* := \min_{p \in \mathbb{R}^m} ~& \frac{1}{2}p^\top B p + \sum_{j = 1}^m p_j \frac{d}{d \tau_j} h_j(\tau_j ) + \sum_{j = 1}^m h_j(\tau_j) \\
\mathrm{s.t.} \quad & a_i(t_i;(\tau_j)_{j=1}^m)^\top p + b_i(t_i;(\tau_j)_{j=1}^m) \leq 0,\nonumber\\
 \quad & \quad \quad \quad \quad\quad \quad\quad \quad ~t_i \in \mathcal{A}_i,~i \in \mathbb{N}_{[1,n_c]}, \\
  %C_z x(t_z^*; (\tau_j)_{j=1}^m)  +  p^\top \nabla_{(\tau_j)_{j=1}^m} \left(C_zx(t_z^*; (\tau_j)_{j=1}^m)\right) \leq c_z, \nonumber\\
%\quad &\forall z \in \{1,\dots,p_x\}, t_z^* \in \argmax_t \left(C_z x(t; (\tau_j)_{j=1}^m) - c_z\right)\\
&p_j \in [\underline{\tau}_j-\tau_j , \overline{\tau}_j-\tau_j ] , \quad j \in \mathbb{N}_{[1,m]},\label{eqn:local_Quad_approxc}\\
&\|p\|_{\infty} \leq \Gamma. %\nonumber%\label{eqn:local_Quad_approxd}
\end{align}
\end{subequations}
The decision variable $p= \begin{bmatrix} p_1&p_2&\cdots&p_m\end{bmatrix}^\top$ is the update direction from the current iterate $(\tau_j)_{j=1}^m$ to a new schedule, $$a_i(t_i;(\tau_j)_{j=1}^m) = \nabla_{(\tau_j)_{j=1}^m} \left(C_i x(t_i; (\tau_j)_{j=1}^m)\right)\in\mathbb{R}^m$$ and $$b_i(t_i;(\tau_j)_{j=1}^m) = C_i x(t_i;(\tau_j)_{j=1}^m) - c_i \in \mathbb{R}$$ are constants, for the given feasible schedule and the finite cardinality sets $\mathcal{A}_i$ contain the maximizers $\mathcal{T}_i^*((\tau_j)_{j=1}^m)$ defined in \eqref{eq:LLP_maximizers} for $i\in\mathbb{N}_{[1,n_c]}$. 
%; i.e., $\mathcal{T}_i^*((\tau_j)_{j=1}^m) \subseteq \mathcal{A}_i[k]$ for each $i \in \mathbb{N}_{[1,n_c]}$. 
The positive scalar $\Gamma$ is a trust region for the current iterate. The matrix $B$ is a symmetric positive definite approximation of the Hessian of the Lagrangian, with multipliers $((\lambda_{il})_{l=1}^{|\mathcal{A}_i|})_{i=1}^{n_c}$, given by
\begin{align} 
& L\left((\tau_j)_{j=1}^m; ((\lambda_{il})_{l=1}^{|\mathcal{A}_i|})_{i=1}^{n_c}\right) \label{eq:Larg} \\
&=
	 \sum_{j=1}^m h_j(\tau_j) + \sum_{i=1}^{n_c} \sum_{l=1}^{|\mathcal{A}_i|} \lambda_{il} \left(C_i x(t_i^{(l)}; (\tau_j)_{j=1}^m) - c_i\right), \nonumber 
\end{align}
where $|\mathcal{A}_i|$ denotes the cardinality of $\mathcal{A}_i$, $i\in\mathbb{N}_{[1,n_c]}$.
% To formulate a quadratic approximation, first consider the Lagrangian of the original problem \eqref{eqn:optim:1} but with constraints only enforced at the maximizers given in \eqref{eq:LLP_maximizers}. This is given by
% where $\lambda_{ik} \in \mathbb{R}$ are the dual variables.
%To understand approximate original problem \eqref{eqn:optim:1} with a quadratic program the lar
%%
% a symmetric positive definite matrix containing information about the curvature of the Lagrangian of the problem.
\begin{remark}
Note \eqref{eqn:local_Quad_approx} is feasible provided the schedule around which the problem is approximated is feasible. This can be seen by setting $p=0$.
\end{remark}

The quadratic program \eqref{eqn:local_Quad_approx} is a local approximation of the original problem \eqref{eqn:optim:1} relative to the schedule $(\tau_j)_{j=1}^m$. The SQP algorithm (SQPA) proceeds by updating this schedule to the schedule $(\tau_j+\gamma p_j^*)_{j=1}^m$, where $p^*$ solves \eqref{eqn:local_Quad_approx} 
% The    
% % To maintain validity of this approximation a trust region is enforced by \eqref{eqn:local_Quad_approxd}. 
% After solving \eqref{eqn:local_Quad_approx} the schedule is updated by stepping in the direction of $p$ by some amount $\gamma$.
and a line search is performed to find the largest $\gamma \in [0,1)$ such that
\begin{equation}\label{eq:SQPFeasCon}
 	g_i^*((\tau_j+\gamma p_j^*)_{j=1}^m) \leq 0, ~ \forall i \in \mathbb{N}_{[1,n_c]},
\end{equation}
and
\begin{equation}\label{eq:LineSearchCon2}
	\sum_{j=1}^m h_j(\tau_j+\gamma p_j^*) \leq \sum_{j=1}^m h_j(\tau_j) + \eta \gamma \sum_{j=1}^m p_j^* \frac{d}{d \tau_j} h(\tau_j),
\end{equation}
where $\eta \in (0,1)$ is a constant algorithm parameter. Corresponding updates of $a_i$, $b_i$, and $\mathcal{A}_i$, $i\in\mathbb{N}_{[1,n_c]}$, $h_j$, $\frac{d}{dt}h_j$, $j\in\mathbb{N}_{[1,m]}$, $B$, and $\Gamma$ are then made to become consistent with the new iterate of the schedule. These steps repeat until a stopping criterion is met. Here
 %for Algorithm~\ref{alg:SQP} 
this criterion relates to the step size becoming smaller than a tolerance $\epsilon^{s}>0$; i.e.,
\begin{equation}
	\gamma \|p^*\|_\infty < \epsilon^{s}. \label{eq:TerminationCond}
\end{equation}

To guarantee a unique solution to \eqref{eqn:local_Quad_approx}, the matrix $B$ must remain positive definite. Since the Lagrangian \eqref{eq:Larg} is non-convex in the schedule a standard Broyden-Fletcher-Goldfarb-Shanno (BFGS) update formula may lead to non-positive definite matrices. Here a damped BFGS update formula from \cite[p.537]{Nocedal2006} is used. The damping ensures a curvature condition is satisfied and that $B$ remains positive definite if the initialization is positive definite. The update formula from \cite[p.537]{Nocedal2006} requires an estimate of the Lagrangian multipliers at current iterate which can be obtained in the process of solving \eqref{eqn:local_Quad_approx}.

The trust region $\Gamma$ is used to reflect confidence in the previous approximation. This can be dynamically updated to allow for larger steps when the approximation is deemed more accurate, and refined as the approximation becomes coarse, in order to reduce the number of iterations needed to perform the line search. Here, \cite[Algorithm 4.1]{Nocedal2006} is used to update the trust region. 
\begin{theorem}\label{the:SQP} 
%Consider Algorithm~\ref{alg:SQP} where 
Given $\epsilon^s>0$ and $\eta\in(0,1)$, suppose the SQPA is initialized with 
$\Gamma>0$, $B$ positive definite, and a schedule that is feasible for \eqref{eqn:optim:1} with cost $f^0$. Then
\begin{enumerate}[label=\it\roman*)]
    \item every schedule iterate is feasible for \eqref{eqn:optim:1}, and
    \item the algorithm terminates in a finite number of steps, such that when the initial solution of \eqref{eqn:local_Quad_approx} is non-zero, 
    \begin{equation}\label{eq:ImprovementIneq}
 	    \sum_{j=1}^m h_j(\tau_j) < f^0
    \end{equation}
    at termination.
\end{enumerate}

% Assume that initial schedule $(\tau_j[0])_{j=1}^m$ with $\tau_j[0]\in \mathcal{D}_j$, $j = \mathbb{N}_{[1,m]}$, satisfies \eqref{eqn:Feasble_Q}, and initial $B_k$ is a positive definite matrix. Then for some finite $\bar{k}$ Algorithm~\ref{alg:SQP}, given  $\epsilon_{step} > 0$ terminates with a schedule $(\tau_j[\bar{k}])_{j=1}^m$ satisfying \eqref{eqn:Feasble_Q} and 

%  with equality holding if and only if $\bar{k}=1$ and $p[0]=0$. Furthermore, if $p[\bar{k}-1] = 0$ and Assumptions~\ref{ass:FiniteMaximizers}, \ref{ass:TiStar} and \ref{ass:Reg} hold the final schedule, $(\tau_j[\bar{k}])_{j=1}^m$, satisfies the first order optimality conditions for \eqref{eqn:optim:1}.
% \begin{enumerate}
% 	\item An initial feasible schedule $(\tau_j[0])_{j=1}^m$ i.e., \eqref{eqn:Feasble_Q} holds and $\tau_j[0]\in \mathcal{D}_j \forall j \in \mathbb{N}_{[1,m]}$, and
% 	\item an initial positive definite matrix $B_0$,
% \end{enumerate}
%  generates a feasible schedule $(\tau_j[k])_{j=1}^m$ that has a nominal objective no worse than the initial schedule i.e.,
%  \begin{equation}
%  	\sum_{j=1}^m h_j(\tau_j[k]) \leq \sum_{j=1}^m h_j(\tau_j[0])
%  \end{equation}.
%  Equality holds if and only if the algorithm terminates in the in the first step with the same initial schedule. Additionally, if algorithm terminates with the optimal step from solving \eqref{eqn:local_Quad_approx} of $p[k] = 0$ and assumptions~\ref{ass:FiniteMaximizers}, \ref{ass:TiStar} and \ref{ass:Reg} hold the final schedule satisfies the first order optimality conditions for \eqref{eqn:optim:1}.
\end{theorem}
\begin{IEEEproof}
See Appendix~\ref{proof:tho:SQP}
\end{IEEEproof}

\begin{remark}
Either the SQPA terminates immediately, with the initial schedule, or a new SIP feasible schedule with strictly improved cost is found in finite steps. By contrast, strict improvement cannot be guaranteed by the penalty methods proposed in \cite{Farokhi2016}. 
\end{remark}

\begin{remark}
If the requested load inputs $v_j$ are restricted to continuous functions, and certain regularity assumptions hold, and the sets $\mathcal{A}_i$ at each step are modified to include points that can be \emph{extended} (in an appropriate sense)
%\footnote{See~\cite{Price1990} for definition of extension.} 
to each of the global maximizers of another schedule that is in the neighborhood of the iterate $(\tau_j)_{j=1}^m$, then it is possible to show that the SQPA, with sufficiently small $\epsilon^{s}>0$, converges to a point that satisfies the first order optimality conditions of \eqref{eqn:optim:1}; i.e., to a local-stationary point. This result is established in~\cite{Price1990}. However, improvement in the cost is the primary aim here, and since it is difficult to show the regularity assumption holds and to ensure $\mathcal{A}_i$ contains the necessary points, the conditions are relaxed to those in Theorem~\ref{the:SQP}, at the expense of not ensuring local optimality. 
\end{remark}

 \begin{remark}
 The algorithm presented in \cite{Price1990} also allows for infeasible initial schedules via an exact penalty method. However, it is observed in practice that this approach is very sensitive to initialization. Within the context of the numerical examples presented in Section~\ref{sec:numerical_results},  initialization of the approach at realistic original flow requests failed to yield a schedule that satisfies \eqref{eqn:Feasble_Q}. This is not unexpected since stationary points for the exact penalty reformulation of the problem are not necessarily feasible. This motivated development of the EHSIPA.
 \end{remark}

\section{Characteristics of the two-stage approach} 
% (fold)
\label{sub:discussion}

%Algorithm~\ref{alg:proposedApp} presents the combined two stage approach to the scheduling problem. The EHSIPA, Algorithm~\ref{alg:EHSIPA}, is first used to generate a feasible schedule which in turn is locally refined by Algorithm~\ref{alg:SQP}. %This leads to the following result.
The two-stage approach to Problem~\ref{prob:scheduling} is as follows. First, the EHSIPA described in Section~\ref{sec:stage1} is applied to generate a schedule that is feasible for \eqref{eqn:optim:1} and $\epsilon^f$-optimal for the restriction associated with a discretization of the decision spaces. In the second stage, this SIP feasible schedule is then improved in cost by application of the SQPA described in Section~\ref{sec:stage2}. 
\begin{theorem} \label{th:MAINRESULT}
Under Assumptions~\ref{ass:feasibleOrig} and~\ref{assum:NLP_Solver}, suppose $\epsilon^{f}, \epsilon^{LBD},$ and $\epsilon^{UBP}$ satisfy \eqref{eq:tol_relations}. 
%Under Assumptions~\ref{ass:feasibleOrig},~\ref{assum:NLP_Solver} then Algorithm~\ref{alg:proposedApp} 
Then the two-stage approach terminates finitely with a schedule $(\tilde{\tau}_j^*)_{j=1}^m$ such that $i)$ $\tilde{\tau}_j^* \in \mathcal{D}_j$, $j \in \mathbb{N}_{[1,m]}$, $ii)$ \eqref{eqn:optim:1:sc} holds, and $iii)$ 
\begin{equation}\label{eq:MainBounds}
	\underline{f} \leq f^* \leq f(\tilde{\tau}_1^*,\dots,\tilde{\tau}_m^*) \leq \hat{f}^* + \epsilon^f,
\end{equation}
where $f^*$ and $\hat{f}^*$ are defined by \eqref{eqn:optim:1} and \eqref{eqn:optim:dd}, respectively, and 
\begin{equation}\label{eq:obLB}
	\underline{f}:= \min_{(\tau_j)_{j=1}^m} f(\tau_1,\dots,\tau_m), \hspace{0.05in} \text{s.t} \hspace{0.05in} \tau_j \in \mathcal{D}_j,~j \in \mathbb{N}_{[1,m]}.
\end{equation}
\end{theorem}
% Assume that $\epsilon^{f}, \epsilon^{LBD},$ and $\epsilon^{UBP}$ are chosen such that \eqref{eq:tol_relations} hold and a continuous algorithm that starts with a feasible schedule, guarantees finite termination with a schedule satisfying \eqref{eq:ImprovementIneq} and \eqref{eqn:Feasble_Q} is used as stage 3 of Algorithm~\ref{alg:proposedApp}. Assume the decision spaces in stage 1 of Algorithm~\ref{alg:proposedApp} are discretized in such a way that assumption~\ref{assum:SIP_point} holds. Then under assumption \ref{assum:NLP_Solver} Algorithm~\ref{alg:proposedApp} used on the scheduling problem defined by Problem~\ref{prob:scheduling} where assumption~\ref{ass:feasibleOrig} holds and $u=u_0$, $u_0 \in \mathbb{R}^{n_u}$
%  terminates finitely with a schedule $(\tau_j^*)_{j=1}^m$ such that $\tau_j^* \in \mathcal{D}_j$, $\forall j \in \{1,\dots, m\}$, \eqref{eqn:Contr} holds and 
% \begin{equation}\label{eq:MainBounds}
% 	\underline{f} \leq f^* \leq f(\tau_1^*,\dots,\tau_m^*) \leq \hat{f}^* + \tilde{\epsilon}^f,
% \end{equation}
% where $f^*$ and $\hat{f}^*$ are defined by \eqref{eqn:optim:1} and \eqref{eqn:optim:dd} respectively.
% \end{theorem}
\begin{IEEEproof}
Proposition~\ref{th:HISPA} gives the right hand side of \eqref{eq:MainBounds}. Theorem~\ref{the:SQP} ensures that the SQPA maintains or improves the cost whilst maintaining SIP feasibility. Hence, the bound and feasibility still hold after stage 2. The left-hand side of \eqref{eq:MainBounds} follows because the cost of any feasible schedule is an upper bound for the optimal cost and that the optimization problem defining $\underline{f}$ is a relaxation of \eqref{eqn:optim:1}; i.e., there are no constraints on the dynamics. 
\end{IEEEproof}

\begin{remark}
Theorem~\ref{th:MAINRESULT} holds for any second stage algorithm 
%if Algorithm~\ref{alg:SQP} used in Algorithm~\ref{alg:proposedApp} was replaced with any algorithm 
that takes an initial schedule that is feasible for \eqref{eqn:optim:1} and terminates after a finite number of steps with a schedule satisfying \eqref{eq:ImprovementIneq} and \eqref{eqn:Feasble_Q}.
\end{remark}
\begin{remark}
The inequality \eqref{eq:MainBounds} bounds the obtained optimality gap. Specifically,  $f(\tau_1^*,\dots,\tau_m^*)-f^* \leq f(\tau_1^*,\dots,\tau_m^*)-\underline{f}\leq \hat{f}^* + \epsilon^f - \underline{f}$.
%\tilde{\epsilon}^f-\underline{f} <  \hat{f}^* +\epsilon^f - \epsilon^{LBP} - \epsilon^{UBP} -\underline{f}$. 
If an upper bound is known for $\hat{f}^*$ then the last expression can be calculated \emph{a priori}. However, in practice these bounds are conservative since $\underline{f} \leq f^*$. This bound relies on calculation of \eqref{eq:obLB} which should be easier than \eqref{eqn:optim:1} since there are no constraints from the dynamics. In particular if $h_j$ are convex then \eqref{eq:obLB} is a convex problem which could be solved using standard convex optimization methods such as those in \cite{boyd2004convex}.
\end{remark}
\begin{remark}
The main result relies on Assumption~\ref{assum:NLP_Solver}. If there exist algorithms that satisfied Assumption~\ref{assum:NLP_Solver} with the original continuous decision spaces, then the HSIPA algorithm for constraint discretization would find the optimal solution, negating the need for the second localized SQPA based refinement stage. The decision spaces are discretized to arrive at tractable problems that satisfy Assumption~\ref{assum:NLP_Solver}. The localized refinement stage can then lead to improvement. 
\end{remark}

Although Theorem~\ref{th:MAINRESULT} guarantees finite termination it does not say anything about the convergence rate or computational complexity and specifically how to choose the initial discretizations of the decision spaces and constraints. The subsequent discussion relates to why a good choice of initial discretizations is not obvious.
\subsection{Effect of Initial Discretizations}
\label{sub:discretizationEffect}
%Algorithm~\ref{alg:proposedApp} 
The two-stage approach based on EHSIPA and SQPA is guaranteed to converge for any choice of initial discretizations of the decision spaces and constraints in \eqref{eqn:optim:1}. However, this choice can have significant impact on both the overall computational time and the achieved optimality gap. 

\subsubsection{Discretization and Computation Time} % (fold)
\label{ssub:discretization_and_computational_time}
The integer programs in the EHSIPA have dimensions linked explicitly to $N_j$ and $T_i$, the cardinalities of the discretized decision spaces and constraint discretization sets. Hence, it may seem intuitive to start with these being small in order to reduce the solve time of these sub-problems. However, sparse initial sets may lead to slower overall computation time. 

To see this, first consider sparse initial constraint discretization sets $\hat{\mathcal{T}}_i$. Using these may lead to an initial lower-bound within the EHSIPA that is far from the optimal, e.g., if the constraint set were empty, then the initial lower bound would be $\underline{f}$. Further, the first call to the upper-bound procedure may require many iterations to add sufficient points to the constraint discretization sets to find a feasible schedule. 

Second, consider sparse initial discretized decision spaces $\hat{\mathcal{D}}_j$. This could yield an infeasible problem \eqref{eqn:optim:dd} and several iterations of the EHSIPA could be required before \eqref{eqn:optim:dd} becomes feasible. Additionally, a sparse decision space may increase the likelihood that the resulting schedule from the first stage is far from a stationary point of the original problem \eqref{eqn:optim:1}. This can increase the convergence time of the SQPA.

These issues are discussed further within the context of a particular example in the numerical results Section~\ref{sec:numerical_results}. 

% subsubsection discretization_and_computational_time (end)
\subsubsection{Discretization and Optimality Gap} % (fold)
\label{ssub:discretization_and_optimality_gap}
The SQPA is guaranteed to improve the initial schedule if possible. However, 
due to non-convexity of \eqref{eqn:optim:1}, the iterates tend towards local stationary points of \eqref{eqn:optim:1}. In general, there is no measure of how far these local stationary points are from the global optimum. Hence, the sub-optimality gap achieved for the two-stage approach is strongly linked to the sub-optimality of the first stage solution; i.e., the distance between $\hat{f}^*$ and $f^*$.
\begin{remark}
Nothing can be said about the achieved optimality gap if the SQPA is initialized at an arbitrary feasible schedule, rather than an optimizer for \eqref{eqn:optim:dd}. 
%It is also hypothesized that optimality of the first stage can result in a superior starting schedule i.e., one that is closer to a local solution with smaller objective, compared to if the first stage just determined any feasible schedule.
\end{remark}
% This also justifies the desire for the first stage of Algorithm~\ref{alg:proposedApp} to find an optimal schedule with respect to \eqref{eqn:optim:dd}, otherwise there would be even less that 

The difference $\hat{f}^*-f^*$ depends on the discretization of the decision spaces. This difference can be made arbitrarily small by adding sufficiently enough points to the discretzized decision spaces $(\hat{\mathcal{D}}_j)_{j=1}^m$. However, this approach is impractical since this may lead to prohibitively large integer programs to be solved as part of the HSIPA. This is highlighted further in Section~\ref{sec:numerical_results}. 

As mentioned in Remark~\ref{rem:GrowthOfDecisionSpaces}, it is of interest to investigate the use of additional information to inform the discretization(s) of smaller decision spaces. Note that it is possible for a smaller discretization to result in a smaller optimality gap. For instance, consider the case where the discretization for each user is chosen with cardinality equal to one, containing only the shift from an optimal schedule for \eqref{eqn:optim:1}; i.e.,  $\hat{\mathcal{D}}_j = \{\tau_j^*\}$ for $j \in \mathbb{N}_{[1,m]}$. In this case, the smallest optimality gap is achieved even though the discretizations are very sparse. 

A possibility in this direction is to re-run the EHSIPA after the first feasible schedule is found, with a finer discretization centered at the solution of the previous iterate such that the cardinality of the decision spaces remain the same.
%, i.e., $N_j[k+1]=N_j[k]$.   
This preliminary idea, along with the effects of the other aforementioned complexities regarding the initial discretizations are demonstrated for a practical example from automated irrigation networks in the subsquent section.

\section{Numerical Results} % (fold)
\label{sec:numerical_results}
In this section, a realistic scheduling setup is investigated numerically based on operational data for an irrigation channel in south-eastern Australia. Results obtained from applying the first stage of the approach developed above are presented for nine instances of the initial discreizations of the decision space and constraints. These initializations range from fine uniform discretizations of both the decision space and the constraints to coarse discretizations of both. The outcomes are compared in terms of feasibility, optimality and computational burden. The results reveal that a dense initialization does not necessarily yield the best outcome. For the example considered it is observed that coarse initialization can yield a feasible schedule that is as good as one obtained from a much denser starting point, without the optimization sub-problems involved becoming overly large. The effectiveness of the second stage is also explored, including comparisons with the penalty based gradient methods from \cite{Farokhi2016}, in terms of improvement in cost and computational burden.

% These different initialization enable the efficacy of the approach to be assessed in a number of ways. These include 
% numerical simulations for an example from automated irrigation networks are presented. The simulations illustrate various properties of the proposed two-stage approach to Problem~\ref{prob:scheduling}. Two primary facets are explored throughout the example. The first feature explored is the how the choice of initial discretizations can effect algorithm performance, highlighting the complexities discussed in Section~\ref{sub:discretizationEffect}. This aspect also highlights how the proposed approach can find feasible schedules with much smaller discretizations and computational times than previous methods that deal with only discrete-time systems and hence uniform sampling. 

% The second feature explored is the suitability of the proposed SQPA for the second, localized refinement, stage of the proposed approach. Specifically, the SQPA is compared, in terms of both optimality and computation time, with modified gradient based algorithms from \cite{Farokhi2016}.
% In this section, numerical examples from irrigation networks are used to illustrate the applicability and some of the properties of the methods presented.
\subsection{Application Gravity Fed Irrigation Networks} % (fold)
\label{sub:application_gravity_fed_irrigation_networks} An irrigation channel is a cascade of pools. Following \cite{WEYER2001, mareels2005systems}, each pool is modeled in the Laplace domain as\vspace{-.04in}
\begin{align*}
y_i(s)=\frac{c_{\mathrm{in},i}}{s} e^{-t_{\mathrm{d},i}s}q_i(s)
-\frac{c_{\mathrm{out},i}}{s}q_{i+1}(s)-\frac{c_{\mathrm{out},i}}{s}o_i(s),
\end{align*}
where $c_{\mathrm{in},i}$ and $c_{\mathrm{out},i}$  (in $1/(\text{min}\sqrt{\text{m}})$) are discharge rates determined by the physical characteristics of the gates used to set the flow between neighbouring pools, and $t_{\mathrm{d},i}$ (in min) is the delay associated with the transport of water along the pool. Here, $o_i(s)$ denotes the overall off-take load on pool $i$, that is, the sum of all the water supplied to the farms connected to this pool. Moreover, $q_i(s)$ represents the head over the upstream gate of pool $i$ raised to the power of $3/2$ which is proportional to the flow (in $\text{m}^3/\text{min}$) of water from pool $i-1$ to pool $i$ and $y_i(s)$ denotes the water level (in m) in pool $i$. For the purpose of this example, the delays are replaced with a first-order Pad\'{e} approximation\footnote{Note that the choice of a first-order Pad\'{e} approximation is justifiable as the pool delays are all parts of closed-loops (with local controllers), with loop-gain cross-overs that are sufficiently small to make the overall closed-loop behavior insensitive to the approximation error~\cite{MichaelCSM}.}. Each pool is controlled, locally, by
$$
q_i(s)=\frac{\kappa_i(\phi_i s+1)}{s(\rho_i s+1)} (u_i(s)-y_i(s)) + \gamma_i q_{i+1}(s),
$$
where $\kappa_i$, $\phi_i$, $\rho_i$ and $\gamma_i$ are appropriately selected control parameters. Furthermore, $u_i(s)$ (in m) denotes the water-level reference signal of pool $i$. The model for a series of pools using the aforementioned model, with first-order Pad\'{e} approximation, can be represented in the form of \eqref{eqn:LTI} where there are four states per pool, two for the controller and two for the level dynamics.

\subsection{Example Parameters and Setup}
\subsubsection{Problem Data}
 Table~\ref{table:poolParams} shows the parameters for the 10 pool system used in this example, which are from validated models of a real channel provided by Rubicon Water Pty Ltd. A fixed $\gamma_i =0.7$ is used for all pools $i \in \mathbb{N}_{[1,10]}$. \al{The reference input $u$ is a set-point to a lower-level feedback controller. In this context it is common that the reference be kept constant as it is desired to maintain a constant level of supply to the offtakes. For this setup} it is fixed to $1m$ for all pools, i.e.,  $u_0 = \mathbbm{1}_{10}$. 
 \begin{table}
\caption{\label{table:poolParams} Parameters, $c_{in,i}$, $c_{out,i}$ (in \textnormal{1/($\text{min} \sqrt{\text{m}})$}), $t_{d,i}$ (in \textnormal{min}), and unitless control parameters $\kappa_i,\phi_i$ and $\rho_i$ for $i\in \mathbb{N}_{[1,10]}$ for the dynamical system used in the simulation.}
\begin{tabular}{|c|c|c|c|c|c|c|}
\hline
Pool No& $c_{\mathrm{in},i}$ & $c_{\mathrm{out},i}$ & $t_{\mathrm{d},i}$ & $\kappa_i$ & $\phi_i$ & $\rho_i$ \\ \hline
$i=1$  & $0.1079$            & $0.1080$             & $1$                & $0.0156$   & $46.648$ & $3.452$ \\ \hline
$i=2$  & $0.0777$            & $0.0777$             & $1.67$             & $0.0091$   & $72.403$ & $5.213$ \\ \hline
$i=3$  & $0.0586$            & $0.0586$             & $2.33$             & $0.0065$   & $99.274$ & $7.084$ \\ \hline
$i=4$  & $0.1269$            & $0.1269$             & $1.67$             & $0.0084$   & $60.305$ & $3.972$ \\ \hline
$i=5$  & $0.0313$            & $0.0313$             & $1.83$             & $0.0092$   & $110.85$ & $8.878$ \\ \hline
$i=6$  & $0.0456$            & $0.0507$             & $4$                & $0.0036$   & $152.73$ & $10.36$ \\ \hline
$i=7$  & $0.0725$            & $0.0725$             & $1.33$             & $0.0119$   & $64.978$ & $4.885$ \\ \hline
$i=8$  & $0.0216$            & $0.0216$             & $3.67$             & $0.0080$   & $147.65$ & $10.28$ \\ \hline
$i=9$  & $0.0366$            & $0.0366$             & $1.67$             & $0.0100$   & $98.231$ & $7.816$ \\ \hline
$i=10$ & $0.2062$            & $0.2331$             & $2$                & $0.0100$   & $48.156$ & $2.101$ \\ \hline
\end{tabular}
%\vspace{-.2in}
\end{table}

The state constraints are envelope constraints on the water level such that
$y_i(t) \in [\underline{y}_i, \overline{y}_i]$ for all $t\in [0,T]$ where $\underline{y}_i=0.9$ (m) and $\overline{y}_i=1.075 (\text{m}) \quad \forall i \in \mathbb{N}_{[1,4]}\cup \{9,10\}$ and $\underline{y}_i=0.88$ (m) and $\overline{y}_i=1.10 \quad \forall i \in  \mathbb{N}_{[5,8]}$.

In this example, there are two users for each pool. Each user $j$, in each pool $i$, has a requested off-take profile  defined by a start time, $s_{ij}$, duration $d_{ij}$ (both in min), and magnitude $m_{ij}$, which is proportional, via a discharge coefficient, to flow (in $\text{m}^3/\text{min}$). This gives typical pulse shape for off-take flows in irrigation networks, i.e., $v_{ij}(t) = m_{ij}(H(t-s_{ij}) - H(t-s_{ij} - d_{ij}))$ where $H:\mathbb{R} \to \{0,1\}$ denotes the continuous time Heaviside step function. Table~\ref{table:distParams} shows the particular off-take load parameters used in this example. The top of Fig.~\ref{fig:UnscheduledVsOptimal_Orders}. displays the requested orders, which follow a realistic demand pattern based on historic orders. It is of note the initial requests cause significant violation of the constraints; see top of Fig.~\ref{fig:UnscheduledVsOptimal_Levels}. 

% The load profiles simulated are pulses defined by a start time, $s_i$, duration, $d_i$, and magnitude, $m_i$, which is a typical shape of a flow load profile within irrigation networks, i.e., $v_i(t) = m_i(H(t-s_i) - H(t-s_i - d_i))$ where $H:\mathbb{R} \to \{0,1\}$ denotes the continuous time Heaviside step function. 

% Two users for each pool are simulated starting at different times, some with different duration and magnitudes, refer to Table~\ref{table:distParams}, where $s_{ij}$ denotes the start time for the demand profile of $j$th user in the $i$th pool, similarly for $m_{ij}$ and $d_{ij}$. This represents a selection of order requests that is consistent with outlet orders for the real pools being modeled. 
The bounds on the allowable shifts are set to $\underline{\tau}_{ij} = -180$ and $\overline{\tau}_{ij} = 180$ for all $i \in \mathbb{N}_{[1,10]}, j \in \{1,2\}$, i.e., the orders can be scheduled by shifting the requested load forward or backwards by up to three hours. The end-user sensitivity is measured with a quadratic function $h_j(\tau) = 0.01 \tau^2$ for all users, i.e., each user experiences a greater cost the further away from the requested start-time. The ideal schedule for each user is to not have their order shifted at all. The planning horizon is set to $T = 1440$min (1 day), for all simulations. \al{In open-channel irrigation networks, it is desirable to reduce the required lead-time for users to request desired off-take profiles. In practice the lead time is typically in the order of days, but there are network operators striving for lead times of around  two hours. Hence, a computational time in the order of an hour for solving the scheduling problem is tolerable.} The aforementioned setup gives all the problem data needed to formulate \eqref{eqn:optim:1}.

\begin{table}\centering
\caption{\label{table:distParams} Parameters for off-take loads used in the simulation. $s_{ij},d_{ij}$ (in \textnormal{min}) and $m_{ij}$ multiplied by a discharge coefficient has units of \textnormal{$\text{m}^3/\text{min}$}}
\begin{tabular}{|c|c|c|c|c|c|c|}
\hline
Pool No & $s_{i1}$ & $d_{i1}$ & $m_{i1}$ & $s_{i2}$ & $d_{i2}$ & $m_{i2}$       \\ \hline
$i=1$   & $200$       & $360$       & $0.0645$    & $500$       & $360$       & $0.0322$ \\ \hline
$i=2$   & $200$       & $180$       & $0.0322$    & $500$       & $180$       & $0.0322$ \\ \hline
$i=3$   & $200$       & $360$       & $0.0322$    & $500$       & $360$       & $0.0645$ \\ \hline
$i=4$   & $200$       & $180$       & $0.0645$    & $500$       & $180$       & $0.0322$ \\ \hline
$i=5$   & $200$       & $360$       & $0.0322$    & $500$       & $360$       & $0.0322$ \\ \hline
$i=6$   & $200$       & $180$       & $0.0290$    & $500$       & $180$       & $0.0580$ \\ \hline
$i=7$   & $200$       & $360$       & $0.0580$    & $500$       & $360$       & $0.0290$ \\ \hline
$i=8$   & $200$       & $180$       & $0.0322$    & $500$       & $180$       & $0.0322$ \\ \hline
$i=9$   & $200$       & $360$       & $0.0322$    & $500$       & $360$       & $0.0645$ \\ \hline
$i=10$  & $800$       & $360$       & $0.0285$    & $500$       & $180$       & $0.0285$ \\ \hline
\end{tabular}
\vspace{-.2in}
\end{table}

% To compare the validity of the LTI model used, a high-fidelity simulation using the St Venant PDE equations, used in industry to accurately simulate these channel systems, is performed both for the unscheduled and scheduled scenario see Fig.~\ref{fig:UnscheduledVsOptimal_LevelsSt}. This shows the LTI model closely matches, with pool 6 showing the most difference, which results in this pool still slightly violating the constraints even with the optimal schedule.

% subsection application_gravity_fed_irrigation_networks (end)

\subsubsection{Algorithm parameters for the EHSIPA and SQPA stages} % (fold)
\label{sub:effect_of_restriction_choice}
In addition to the parameters mentioned throughout Section~\ref{sec:stage1}, the HSIPA requires specification of the tolerance for problem \eqref{eq:RES}, $\epsilon^{RES}\in \mathbb{R}_{>0}$, an initial tolerance for \eqref{eq:LLP}, $\epsilon^{LLP}\in \mathbb{R}_{>0}$, which is refined throughout the HSIPA by steps $\epsilon^{LLP} \gets \epsilon^{LLP}/r^{LLP}$ for specified $r^{LLP}>1$, as detailed in \cite[Algorithm 2]{Djelassi2017}. Each parameter is independently varied to determine a set of parameters that give reasonable and reliable computational performance. This results in the set of parameters chosen as
$(\epsilon^f,\epsilon^g,r^g,r^{LLP},\epsilon^{LBP},\epsilon^{UBP},\epsilon^{RES},\epsilon^{LLP},l_{max},k_r^{max}) = (5,10^{-3},1.5,1.2,0.5,0.5,10^{-4},0.01,10,10)$. The focus of the simulation study is directed to the effect of the initial discretization of both the decision spaces and constraints.  

The parameters for the SQPA are chosen similarly, resulting in parameters $(\Gamma,\eta, B, \epsilon_{step}) = (5,0.33,I_m,0.5)$. 

\subsubsection{EHSIPA Initializations}
The discretizations $\hat{\mathcal{D}}_{ij}$ are chosen to be a uniformly sampled version of the continuous sets $\mathcal{D}_{ij}=[-180,180]$ with sample period $\Delta_\tau$, i.e., $\hat{\mathcal{D}}_{ij} = \{\underline{\tau}_{ij},\underline{\tau}_{ij}+\Delta_\tau,\dots,\overline{\tau}_{ij}\}$ for all $i \in \mathbb{N}_{[1,10]}, j \in \{1,2\}$, and the initial discrete constraint sets $\hat{\mathcal{T}}_i$ are chosen to be uniformly sampled time points with a period $\Delta$. Nine different simulation study configurations are considered as outlined below:
\begin{enumerate}
	\item $(\Delta_\tau,\Delta) = (60,15)$;
	\item $(\Delta_\tau,\Delta) = (30,15)$;
	\item $(\Delta_\tau,\Delta) = (15,15)$;
	\item $(\Delta_\tau,\Delta) = (5,15)$;
	\item $(\Delta_\tau,\Delta) = (15,\infty)$ (i.e., initial constraint set is empty);
	\item $(\Delta_\tau,\Delta) = (15,30)$;
	\item $(\Delta_\tau,\Delta) = (15,1)$;
	\item First the EHSIPA is run with $(\Delta_\tau,\Delta) = (30,15)$ then the EHSIPA is rerun with discretization centered around the optimal solution with $(\Delta_\tau,\Delta) = (2.5,15)$. Specifically, for user $j$ on pool $i$ denote $\tau_{ij}^*$ as the optimal shift from configuration 2. The restricted set is chosen as $\hat{\mathcal{D}}_{ij} = \{\tau_{ij}^*-30,\tau_{ij}^*-27.5,\dots,\tau_{ij}^*,\dots,\tau_{ij}^*+30\}$;
	\item $(\Delta_\tau,\Delta)=(15,0.25)$ with a single iteration of \eqref{eq:LBD} only.
\end{enumerate}
The first four configurations are used to explore the effect of the decision space discretizations, while configurations 3 and 5-7 explore the effect of the initial constraint discretization. Configuration 8 examines the potential for an improved method of decision space updates as discussed in Section~\ref{ssub:discretization_and_optimality_gap}. Configuration 9, which does not involve update of the initial discretizations, is equivalent to the discrete-time system methods in \cite{5400193}. It is included for comparison with the approach proposed here.

% 	\item $(\Delta_\tau,\Delta)=(15,0.25)$ but only single iteration of \eqref{eq:LBD}.
% For scenario 9) the schedule is computed by solving \eqref{eq:LBD} with uniform sampled contraints with sample period of $0.25$.
% Algorithm~\ref{alg:proposedApp} requires as an input an initial discretization for $\hat{\mathcal{D}}_{ij}$ for each user $j$ and pool $i$ and an initial sampling of the constraints $\hat{\mathcal{T}}_i$. Each of these variables are varied independently whilst keeping the other constant to explore their effect on computational runtime and result for both stages of the Algorithm~\ref{alg:proposedApp}. Nominally, $\hat{\mathcal{D}}_{ij}$ is chosen to be uniformly sampled with sample period $\Delta_\tau$ i.e., $\hat{\mathcal{D}}_{ij} = \{\underline{\tau}_{ij},\underline{\tau}_{ij}+\Delta_\tau,\dots\overline{\tau}_{ij}\}$ for all $i \in \{1,\dots,10\}, j \in \{1,2\}$ and the initial discrete sets $\hat{\mathcal{T}}_i$ are chosen to be uniformly sampled time points with a period $\Delta$. The values of $\Delta_\tau$ and $\Delta$ explored are given in Table~\ref{table:Stage1_Params}.
% \begin{table}
% \caption{\label{table:Stage1_Params} Numerical parameters for the stage 1.}
% \centering
% \begin{tabular}{|c|c|c|}
% \hline
% Parameter & Values trialled & Nominal Value \\\hline
% $\Delta_\tau$ & $\{60,30,15,5,2.5^*\}$ & 15 \\\hline
% $\Delta$ & $\{\infty,30,15,1\}$ & 15 \\\hline
% \end{tabular}
% \end{table}

\subsection{Implementation of proposed two-stage approach} % (fold)
\label{ssub:stage_2_of_algorithm_alg:proposedapp}
\subsubsection{EHSIPA implementation} % (fold)
\label{ssub:extended_hsipa_stage_algorithm_alg:ehsipa}
% subsubsection extended_hsipa_stage_algorithm_alg:ehsipa (end)
\al{Numerical solver \texttt{ode45} in Matlab is used to evaluate the integrals in the components $\bar{x}_{ij}^v(t)$, $\bar{x}_0(t)$,  $i \in \mathbb{N}_{[1,10]}, j\in \mathbb{N}_{[1,2]}$ required to compute \eqref{eq:LLP}.}
% To evaluate \eqref{eq:LLP}, the components $\bar{x}_{ij}^v(t)$, $\bar{x}_0(t)$,  $i \in \mathbb{N}_{[1,10]}, j\in \mathbb{N}_{[1,2]}$, are calculated for each $t\in \{0,\delta_s,2 \delta_s,\dots,T\}$ using the \texttt{ode45} function in Matlab with a specified sample period, $\delta_s$.
By exploiting linearity and time-invariance of \eqref{eqn:LTI}, each constraint $k$ is calculated for a given schedule by summing $\bar{x}_0(t)$ with shifted versions of $\bar{x}_j^v(t)$ and multiplying the result by the corresponding row $C_k$. The integrals are evaluated on a uniform discretized set with an initial sample period chosen as $\delta_s=0.1$. If, during the HSIPA, the maximum difference between any two consecutive samples of the constraint is greater than the current required tolerance $\epsilon^{LLP}$, the sample period $\delta_s$ is reduced accordingly and the intergrals are re-evaluated on discretization set with the smaller sample period. 

The sub-problems in the HSIPA are formed as the linear (mixed-)integer programs described in Section~\ref{ssub:mixed_integer_linear_program} and solved using Gurobi~\cite{gurobi}, interfaced with Matlab.

  \subsubsection{Localized refinement stage} % (fold)
  \label{ssub:localized_refinement_stage}
  For the second localized refinement stage the following three algorithms are compared: 
  \begin{enumerate}[label=\roman*)]
      \item the proposed SQPA, described in Section~\ref{sec:stage2};
      \item the projected gradient algorithm in \cite{Farokhi2016} with log-penalty and $\epsilon=1500$;
      \item the projected gradient algorithm in \cite{Farokhi2016} with exponential-penalty and $\vartheta=1500$.
  \end{enumerate}
  The line search in the projected gradient algorithm in \cite{Farokhi2016} used for ii) and iii) is modified slightly to include a feasibility condition similar to \eqref{eq:SQPFeasCon}. This allows the comparison to focus on the final objective value and computation time. 
  The derivatives required to implement the SQPA can be constructed via Lemma~\ref{lem:constraintDef} from simulations of the linear system dynamics, as used to evaluate \eqref{eq:LLP} in the EHSIPA.
  By contrast, the derivatives of the augmented penalty terms needed for the penalty based algorithms are computed using a first order backward finite difference method; i.e.,
  %as follows: 
  the partial derivative of given function $l(\cdot)$ with respect to shift $\tau_{\ell}$ at the current schedule $(\tau_j)_{j=1}^m$ is approximated as
  \begin{equation}
	\frac{\partial l}{\partial \tau_\ell} \approx \frac{l((\tau_j)_{j=1}^m) - l((\tau_j)_{j\neq \ell};\tau_{\ell} - \delta_s)}{\delta_s},
\end{equation}
with $\delta_s$ set to be smaller than the final sample period from all simulations required for the EHSIPA (for this example $\delta_s=0.01$). This approximation is computationally much more efficient and numerically robust, compared to calculating the derivatives of the penalty functions on the basis of explicit formulae, as these do not directly correspond to a simulation of the linear dynamics.
%involve the continuous integrals and matrix exponentials in the explicit terms \eqref{eqn:barx0} and \eqref{eqn:barxv}.

All simulations were carried out with Matlab 2018b on a Windows PC with 16GB RAM and Intel Core i7-4790K CPU @4.00GHz processor.
\subsection{Results} % (fold)
\label{sub:results}
\begin{figure}\centering
\includegraphics[width=1.0\linewidth]{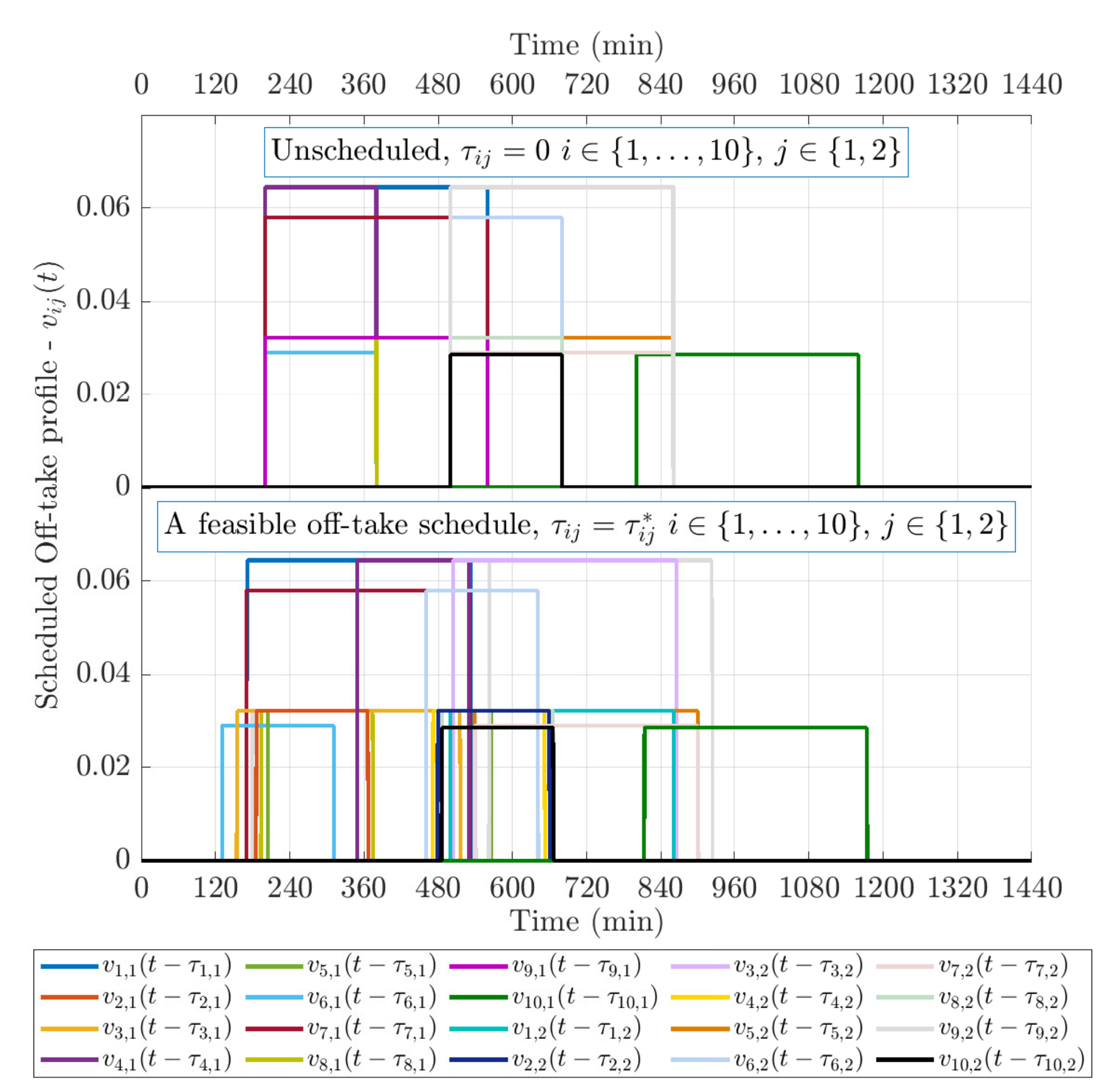}
\caption{\label{fig:UnscheduledVsOptimal_Orders} Top: Unscheduled 20 off-take profiles. Bottom: Off-take profiles under a feasible schedule using proposed approach with configuration 4 followed by SQPA. This feasible sub-optimal schedule is close to original, as per users desire.}
\vspace{-.1in}
\end{figure}
\begin{figure*}\centering
\includegraphics[width=0.78\linewidth]{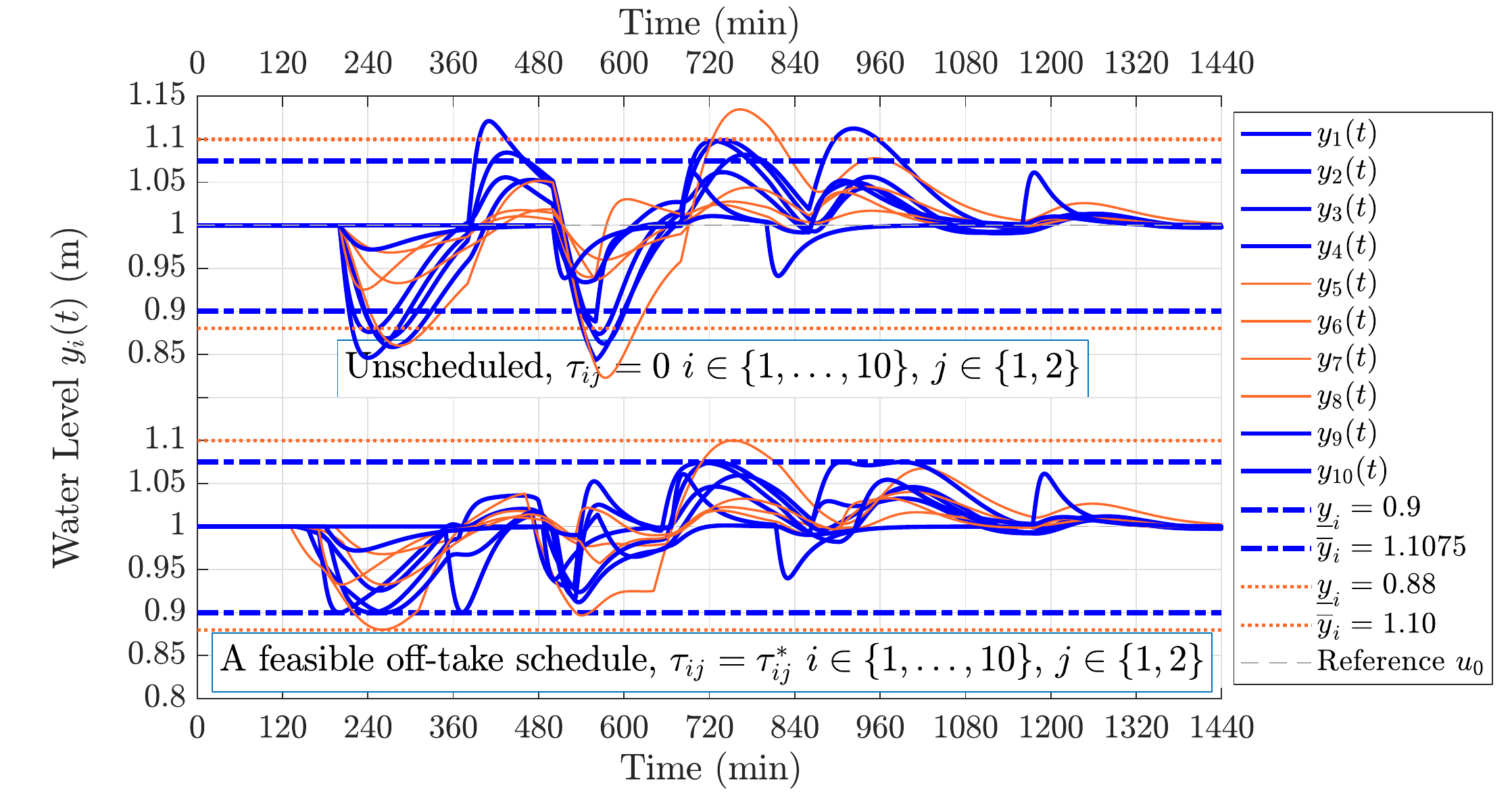}
\caption{\label{fig:UnscheduledVsOptimal_Levels} Levels $y_i(t)$ (in m) in response to nominal off-takes, where constraints are clearly violated on both the upper and lower bounds (top figure); and feasible schedule from configuration 4, where the levels are within the constraint limits (bottom figure).}
\vspace{-.2in}
\end{figure*}
\begin{figure}\centering
\includegraphics[width=0.95\linewidth]{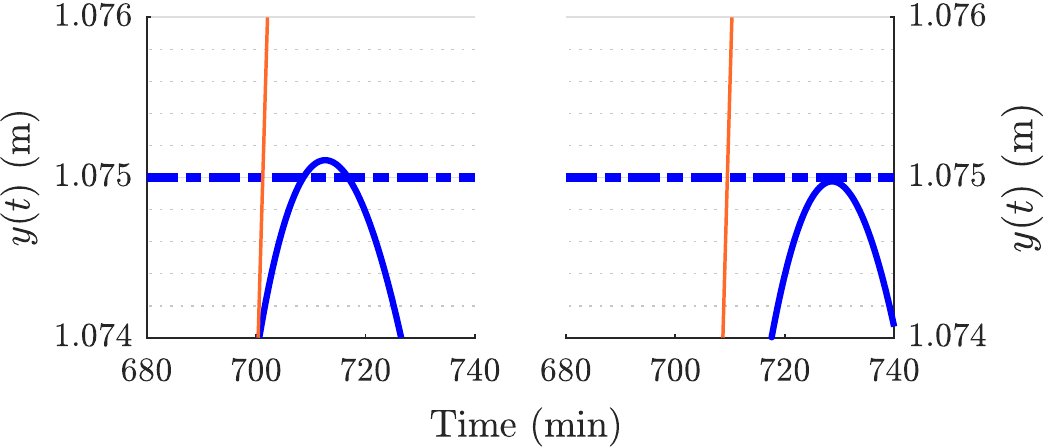}
\caption{\label{fig:constraintSatisfaction} Levels $y_i(t)$ in response to the off-take schedule using the uniform discretization method from \cite{5400193}, i.e., configuration 9, (left), which does not meet constraints, and a off-tale schedule from configuration 3 combined with SQPA, where the levels are within the constraint limits.}
\vspace{-.2in}
\end{figure}

\begin{figure}
    \centering
\includegraphics[width=0.85\linewidth]{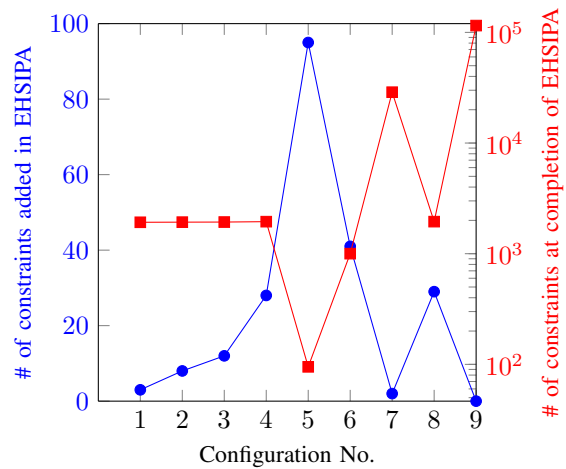}
    \caption{Number of constraints added throughout the EHSIPA algorithm and the total number of constraints for the 9 different configurations.}
    \label{fig:constraints}
    \vspace{-0.1in}
\end{figure}

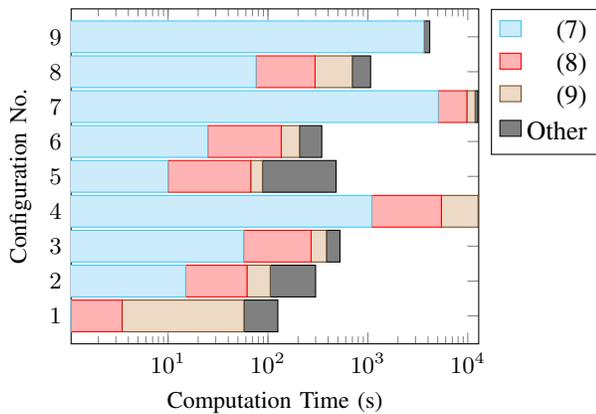
\begin{figure}
    \centering
\begin{tikzpicture}
\begin{semilogxaxis}[
    xbar stacked,
	bar width=12pt,
	legend cell align=right,
legend pos=outer north east,
	%nodes near coords,
    enlarge x limits=false,%0.15,
    %legend style={at={(500,1.0)},
     % anchor=north,legend columns=-1},
    xlabel={\small Computation Time (s)},
    ylabel={\small Configuration No.},
    tick label style={font=\small}, 
    ytick=data,
    %x tick label style={rotate=45,anchor=east},
    ]
\addplot plot coordinates {(1.07,1) (15,2) 
  (57,3) (1095,4) (10,5) (25,6) (5079,7) (76.1,8) (3641.9,9)};
\addplot plot coordinates {(2.4,1) (46.6,2) 
 (213.1,3) (4355,4) (57.2,5) (110.9,6) (4707,7) (219.6,8) (0.001,9)};
\addplot plot coordinates {(54,1) (43.7,2)
  (114.4,3) (7234,4) (20.9,5) (70.9,6) (2041,7) (402,8) (0.001,9)};
\addplot plot coordinates {(67.93,1) (194.2,2) 
  (140.7,3) (170,4) (392.2,5) (138.4,6) (797,7) (366,8) (515.1,9)};
\legend{\eqref{eq:LBD}, \eqref{eq:UBD}, \eqref{eq:RES}, \strut Other}
\end{semilogxaxis}
\end{tikzpicture}
    \caption{The breakdown of the total runtime for the EHSIPA for each configuration.}
    \label{fig:runtimes}
    \vspace{-0.1in}
\end{figure}

The top of Fig~\ref{fig:UnscheduledVsOptimal_Orders} shows the requested (unscheduled) off-take profiles and the bottom shows the off-take profiles shifted with a feasible sub-optimal schedule obtained from the proposed algorithm \al{initialized according to} configuration 4. The corresponding water level trajectories are displayed in Fig.~\ref{fig:UnscheduledVsOptimal_Levels}. The shifted \al{off-takes} remain ``close'' to the requested ones under the sub-optimal schedule, and the levels remain within the constraints. Achieving constraint satisfaction relies on the first stage, EHSIPA, to find a feasible schedule. From Fig.~\ref{fig:constraintSatisfaction} it can be seen that the purely discrete-time (uniform sampling) approach from \cite{5400193} is not able to meet the constraints even with fine discretization (i.e., $\Delta=0.25$, configuration 9). \al{By contrast}, the EHSIPA is able to yield a feasible solution with only 95 constraints (configuration 5); see Fig.~\ref{fig:constraints}. In addition, \al{in} configurations 3, 5 and 6, which have the same \al{initial} decision spaces, \al{the proposed algorithm is} almost a whole order of magnitude faster than \al{in} configuration 9; see Fig.~\ref{fig:runtimes}. 

\begin{figure}\centering
\includegraphics[width=0.85\linewidth]{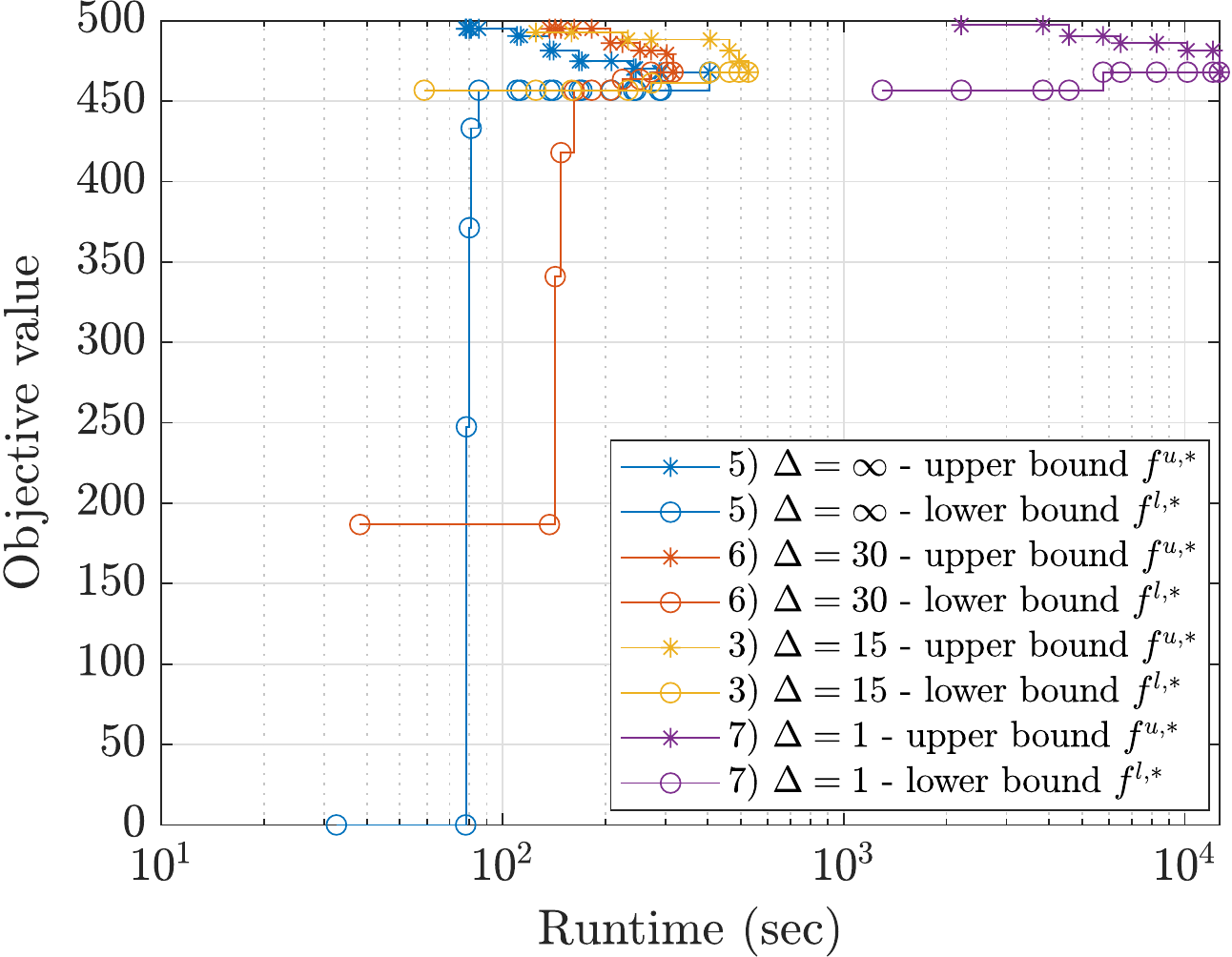}
\caption{\label{fig:initCon} Figure showing convergence of upper and lower bounds of EHSIPA, for different initial discretizations of $\hat{\mathcal{T}}_i$. A less dense discretization, $\Delta = \infty$, results in more a conservative initial lower bound. Whereas, a dense discretization, $\Delta = 1$ takes much longer to calculate initial lower bound.}
\vspace{-.2in}
\end{figure}

The runtime of the EHSIPA increases in proportion to the size of the decision spaces as highlighted by configurations 1-4 in Fig.~\ref{fig:runtimes}. \al{Configurations 3 and 5-7, have $N_j=24$ possible choices for each of the 20 off-takes, which corresponds to $24^{20}$ possible combinations. As such, it is intractable to find the optimal solution via a brute-force (exhaustive search) approach. This highlights the importance of formulating the sub-problems as (mixed-) linear-integer programs as such problems are amendable to widely available dedicated solvers.} Configurations with the same initial decision spaces, i.e., 3 and 5-7, are used to explore the effect of the initial constraint discretization. Of these, configuration 6, with $\Delta=30$, resulted in the fastest run time of the first stage (EHSIPA).  Configuration 5 corresponds to the coarsest initial and final discretization, but more constraints are added during the algorithm execution. As such, it requires many more calls to the lower-level problem \eqref{eq:LLP} and iterations of the sub-problems. Note in particular, that the initial lower-bound is much worse for configurations 5 and 6, as shown in Fig.~\ref{fig:initCon}. This accentuates the points discussed in Section~\ref{ssub:discretization_and_computational_time}.

\begin{figure*}
 \centering
\includegraphics[width=0.85\linewidth]{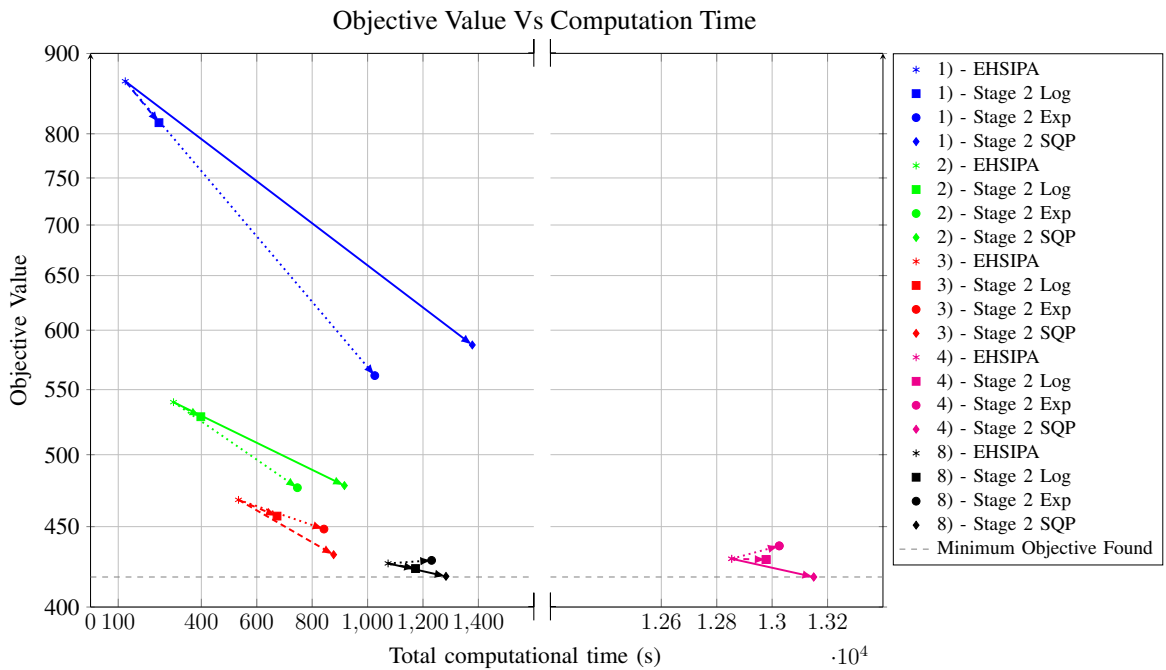}
\vspace{-0.05in}
 \caption{\label{fig:runtimeVsObj} Runtime vs objective for the different choices of $\hat{\mathcal{D}}_{ij}$ and different algorithms for stage 2 of proposed approach}
 \vspace{-.2in}
 \end{figure*}

The configurations with finer discretization resulted in smaller objective $\hat{f}^*$ obtained from EHSIPA; see Fig.~\ref{fig:runtimeVsObj}. Configuration 8 results in a final objective that is $0.67\%$ lower than configuration 4, however it achieves this in $8.4\%$ of the time. This highlights the potential for further improvements to decision space update methods.

The arrows in Fig.~\ref{fig:runtimeVsObj} show the impact of the second stage on the overall objective and computation time. Both a longer run time and larger improvement in objective is made with stage 2 for configuration 1. The logarithmic penalty method was the fastest for all configurations but only resulted in marginal improvement of the objective in each case. The proposed method resulted in improvement in all configurations and had the greatest improvement for configurations 3, 4 and 8. Interestingly, for configuration 4 and 8, the exponential penalty method resulted in an increase in the overall objective, shown by the positive gradient of the dotted arrow.  That is, the penalty method cannot guarantee a strict improvement of the schedule, unlike the SQPA method; see \eqref{eq:ImprovementIneq}.

\begin{figure*}\centering
\includegraphics[width=0.78\linewidth]{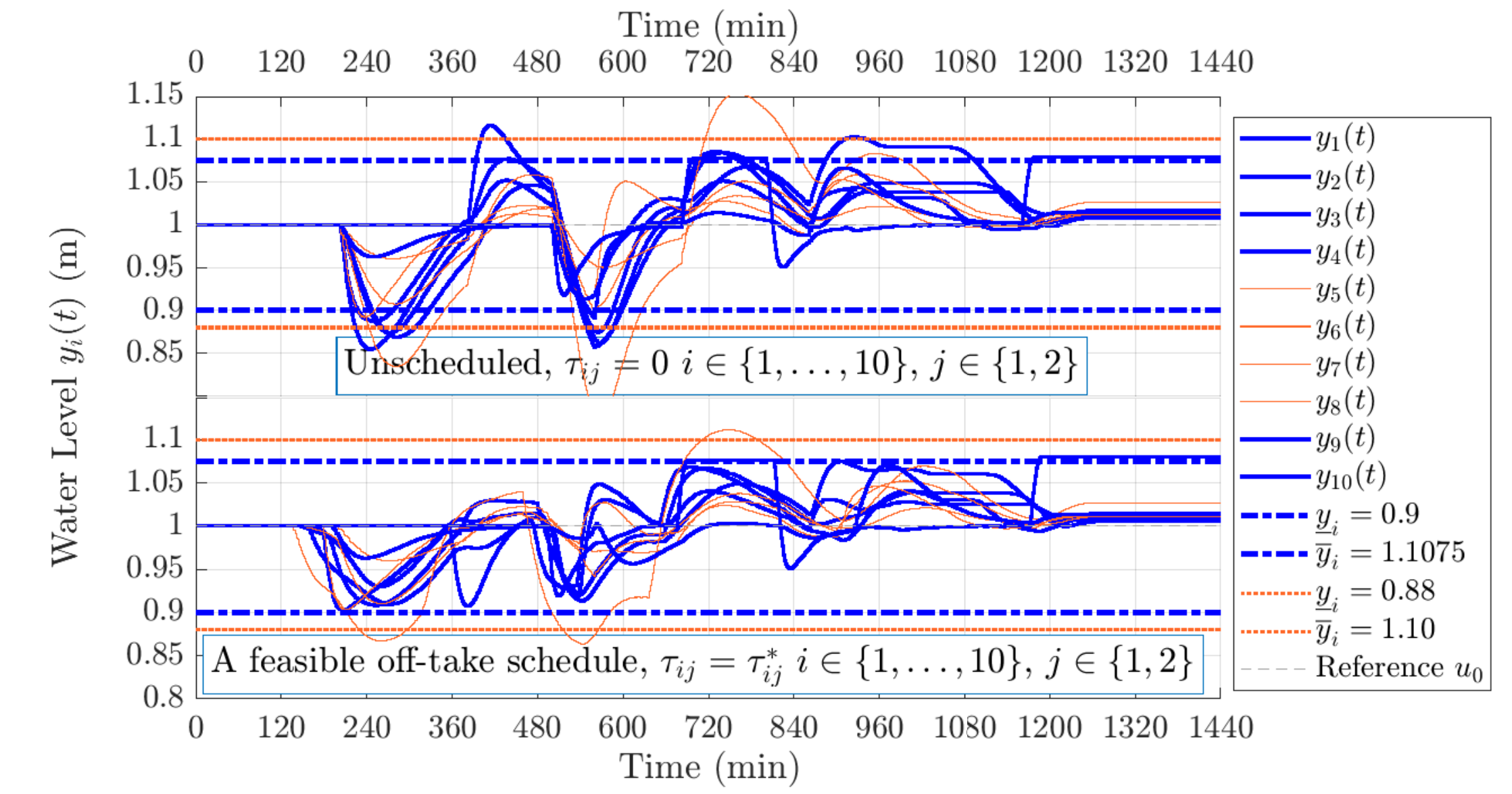}
\vspace{-0.05in}
\caption{\label{fig:UnscheduledVsOptimal_LevelsSt} High fidelity St Venant PDE models based simulation of levels $y_i(t)$ in response to nominal off-takes, where constraints are clearly violated on both the upper and lower bounds (top figure); and a feasible off-take schedule from configuration 4, where the levels are almost all within the constraint limits (bottom figure).}
\vspace{-.1in}
\end{figure*}

Finally, to highlight the applicability of the proposed method in practice, the water levels for the unscheduled and a feasible sub-optimal scheduled off-take loads are simulated using high fidelity models which include, low-level control system actuator saturation and the St Venant PDE model for the open-water dynamics \cite{chaudhry2007open}. Fig.~\ref{fig:UnscheduledVsOptimal_LevelsSt}. shows the optimal schedule only slightly violates the constraints for only two of the pools.


\section{Conclusions and Future Work} \label{sec:conclusion}
 A two-stage approach to finding a sub-optimal feasible solution to a rigid-profile input scheduling problem for a continuous-time linear time-invariant dynamical system is proposed. The first stage pertains to the discretization of both the decision spaces and constraints. Through iteratively refining the discretizations and solving tractable (mixed-)integer linear programs as sub-problems, a schedule that is feasible for the original continuous-time problem is generated. This schedule is then improved by the second stage in a way that guarantees improvement in the objective and feasibility with respect to the original problem. A sequential quadratic programming method is applied for this stage and shown to meet both of these requirements.  Simulation results from a realistic example from automated irrigation networks illustrate the advantages and some of the properties, trade-offs for the proposed algorithm(s). Future work is to explore how prior information or multiple iterations of the first stage of the algorithm can be used to inform more suitable choices of restricted decision spaces. Other work can be done on characterizing the convergence for the various algorithms, and development of feasible methods for computation of useful optimality gaps to detect acceptable solutions. Finally, further work could explore how to update the schedule, perhaps on a receding horizon basis, to accommodate for model mis-matches and to also account for changes to requested orders. 
 
%  to determine an initial  uses discretization to formulate . It iteratively refines the discretizations
%  This paper considers load scheduling and demand management on a continuous-time linear time-invariant dynamical system. This problem is formulated as a semi-infinite program. A two-stage algorithm is proposed to find sub-optimal feasible schedules. The first stage requires an initial discretization of decision sets and constraint sets. It then successively refines these discretizations by extending a hybrid-SIP algorithm which solves a sequence of sub-problems, which are formulated as (mixed-)integer linear programs, to arrive at a feasible schedule that is optimal over the sampled decision space in finite time. The final stage of the algorithm uses the schedule from the first stage and employs an algorithm to locally refine this schedule whilst maintaining feasibility. A sequential quadratic programming algorithm was applied for this stage and shown to yield improvement in the schedule whilst maintaining feasibility. A realistic example from irrigation networks was used to illustrate, compare and discuss some of the properties, trade-offs for the algorithm(s). Future work is to explore how to use prior information or multiple iterations of the first stage of the algorithm to inform more suitable choices of restricted decision spaces. Other work can be done in characterizing the convergence for the various algorithms, and development of feasible methods for computation of useful optimality gaps to detect acceptable solutions.
 \bibliographystyle{ieeetr}
 \bibliography{RigidDemandScheduling}

 \appendices
\section{Proofs of Lemmas and Theorems}
\subsection{Proof of Lemma~\ref{lemm:RESLB}} % (fold)
\label{sec:proof_of_lemma_lemm:reslb}
\al{Assume $\eta^*<0$. Since, the objective of \eqref{eq:RES} is $-\eta$ then optimality of $\eta^*<0$ implies that} for $\eta=0$ there are no feasible schedules $(\tau_j)_{j=1}^m$ for \eqref{eq:RES}. Let $\mathcal{R}:= \{(\tau_j)_{j=1}^m : \tau_j \in  \hat{\mathcal{D}}_j, \,\forall j\in\mathbb{N}_{[1,m]}, \sum_{j=1}^m h_j(\tau_j) \leq f^{R}\}$, which is non-empty since $f^R$ is chosen to be greater than a lower-bound $\hat{f}^{*,L}$ as per \cite[Algorithm 2]{Djelassi2017}. Then, \al{since $\eta=0$ is not feasible for \eqref{eq:RES}}, for all $(\tau_j)_{j=1}^m \in \mathcal{R}$ there exists a constraint $i \in \mathbb{N}_{[1,n_c]}$ and $t \in \hat{\mathcal{T}}_i$ such that
\begin{equation}\label{eq:RESProof}
	C_ix(t; (\tau_j)_{j=1}^m) - c_i > 0. 
\end{equation}
 Therefore, since $\hat{\mathcal{T}}_i\subset \mathcal{T}_i$, there is no schedule in $\mathcal{R}$ that satisfies \eqref{eqn:optim:dd:sc}. Hence $f^{R} < \hat{f}^*$ which completes the proof.
% The proof is constructed by showing that all schedules in the discretized decision space satisfying \eqref{eq:RESfREScon} are infeasible when the condition $\eta^*<0$ holds. Let all schedules in the discretized decision space satisfying \eqref{eq:RESfREScon} be defined by $\mathcal{R}:= \{(\tau_j)_{j=1}^m : \tau_j \in  \hat{\mathcal{D}}_j, \,\forall j\in\mathbb{N}_{[1,m]}, \sum_{j=1}^m h_j(\tau_j) \leq f^{RES}\}$. By assumption $\eta^*<0$ hence, \eqref{eq:RES} must be infeasible for $\eta=0$, otherwise there is the contradiction since zero is more optimal than $\eta^*$, i.e., $\eta=0>\eta^*$. For this to hold then every schedule $(\tau_j)_{j=1}^m \in \mathcal{R}$ must 	not satisfy \eqref{eq:RESCicon} with $\eta=0$ for at least a constraint $i \in \mathbb{N}_{[1,n_c]}$ and $t \in \hat{\mathcal{T}}_i$, owing to these schedules by definition satisfying \eqref{eq:RESfREScon}.
% This can be stated as for all $(\tau_j)_{j=1}^m \in \mathcal{R}$ there exists a constraint $i \in \mathbb{N}_{[1,n_c]}$ and $t \in \hat{\mathcal{T}}_i$ such that
% \begin{equation}\label{eq:RESProof}
% 	C_ix(t; (\tau_j)_{j=1}^m) - c_i > 0. 
% \end{equation}
% Hence, there exists a constraint $i \in \mathbb{N}_{[1,n_c]}$ such that $g_i^*((\tau_j)_{j=1}^m)) > 0$ for every schedule $(\tau_j)_{j=1}^m \in \mathcal{R}$ and therefore there are no schedules in $\mathcal{R}$ that are SIP feasible, i.e., no schedules satisfying \eqref{eq:RESfREScon} and \eqref{eqn:Feasble_Q}. Hence, $f^{RES} < \hat{f}^*$ which completes the proof.

% section proof_of_lemma_lemm:reslb (end)
\subsection{Proof of Lemma~\ref{lem:decisionHalve}} % (fold)
\label{sec:proof_of_lemma_lem:decisionhalve}
Let $(\tau_j^{c})_{j=1}^m$ be a schedule satisfying Assumption~\ref{ass:feasibleOrig}; i.e., a strictly feasible schedule. Then there exists $\gamma>0$ such that  
\begin{equation}
	g_i^*((\tau_j^{c})_{j=1}^m) \leq -\gamma.
\end{equation}
In view \eqref{eqn:barx0} and \eqref{eqn:barxv}, the dependence of the left-hand side of the constraint \eqref{eqn:optim:1:sc} is continuous in $(\tau_j)_{j=1}^m$, whereby $g_i^*$ is continuous as the max of continuous functions. Therefore, for every $\epsilon>0$ there exists $\delta(\epsilon) > 0$ such that 
\begin{subequations}
\begin{align}
	&\max_j \|\tau_j - \tau_j^{c}\| < \delta(\epsilon) \label{eq:lemConConstrainta}\\
	&\qquad \Rightarrow g_i^*((\tau_j)_{j=1}^m) < g_i^*((\tau_j^{c})_{j=1}^m) + \epsilon
	\leq -\gamma + \epsilon.\label{eq:lemConConstraintc}
\end{align}
\end{subequations}
In particular, with $\epsilon = \gamma$, if given schedule $(\tau_j)_{j=1}^m$ satisfies the corresponding condition \eqref{eq:lemConConstrainta}, then it is feasible for the original problem by \eqref{eq:lemConConstraintc}. 

Define the distance between $\tau_j^{c}$ and given $\hat{\mathcal{D}}_j$ as $d_j = \min_{\tau_j \in \hat{\mathcal{D}}_j} |\tau_j - \tau_j^{c}|$, and let $\bar{d} = \max_j d_j$. It follows by construction of the update \eqref{eq:disDUpdate} that this distance is halved with each update. To achieve feasibility the distance must be made smaller than $\delta(\gamma)$, which from an initial value of $\bar{d}_0$ takes $\ceil{\log_2\frac{\bar{d}_0}{\delta(\gamma)}}$ updates.  
%If at iterate $k$, \eqref{eqn:optim:dd} is feasible on decision spaces $(\hat{\mathcal{D}}_j[k])_{j=1}^m$ then proof is complete. Otherwise, the discretization is updated per \eqref{eq:disDUpdate}. Thus, $d_j[k+1] \leq \frac{d_j[k]}{2}$ for all $j\in \mathbb{N}_{[1,m]}$ and hence $\overline{d}[k+1]\leq \frac{\overline{d}[k]}{2}$. Hence, if $\overline{k} = \ceil{\log_{2}\frac{\overline{d}[0]}{\delta_{\gamma}}}$ then $\overline{d}[\overline{k}]< \delta_{\gamma}$. This concludes the proof.

% Let $(\tau_j^{(k)})_{j=1}^m$ be the schedule that is minimum distance from $(\tau_j^{c,*})_{j=1}^m$ i.e, $|\tau_j^{(k)}-\tau_j^{c,*}| = d_j^{(k)}$ for $i\in\mathbb{N}_{[1,m]}$. Then $\max_j \|\tau_j^{(k)} - \tau_j^{c,*}\| = \overline{d}^{(k)} < \delta_{\gamma}$ for $\overline{k} = \ceil{\log_{2}\frac{\overline{d}^{(0)}}{\delta_{\gamma}}}$ hence $(\tau_j^{(\overline{k})})_{j=1}^m$ satisfies \eqref{eq:lemConConstrainta} with $\epsilon=\gamma$ and hence is a feasible schedule. 

\subsection{Proof of Lemma~\ref{lem:constraintDef}} % (fold)
\label{proof:lem:contraintDef}
Note that
\begin{align*}
\frac{\partial}{\partial \tau_\ell} &\left(C_i \bar{x}_0(t)+C_i\sum_{j=1}^m\bar{x}_j^v(t-\tau_j)\right)\\
=&\frac{\partial}{\partial \tau_\ell} C_i \int_0^{t - \tau_\ell} \exp(A(t- \tau_\ell - \beta ))E_\ell v_\ell (\beta) d \beta\\
=&- C_i E_\ell v_\ell(t - \tau_\ell) \\
&+ C_i \int_0^{t - \tau_\ell} \frac{\partial}{\partial \tau_\ell} \exp(A(t- \tau_\ell - \beta ))E_\ell v_\ell (\beta) d \beta\\
=&- C_i E_\ell v_\ell(t - \tau_\ell)  \\
& - C_i \int_0^{t- \tau_\ell} A \exp (A (t - \tau_\ell - \beta)) E_\ell v_\ell(\beta) d \beta\\
\end{align*}
Proof follows from piecewise continuity of $v_j(t - \tau)$.

% section proof_of_lemma_ (end)
\subsection{Proof of Theorem~\ref{the:SQP}}
\label{proof:tho:SQP}
% \subsubsection{First statement - feasibility}
To facilitate the development, relevant variables are indexed by iteration $k$ in what follows. The initial iteration state is $k=0$. Each iteration of the SQPA increments the index. 

{\it i)} 
Assume that $(\tau_j[k])_{j=1}^m$ satisfies \eqref{eqn:Feasble_Q} and $\tau_j[k] \in \mathcal{D}_j$ for all $j \in \mathbb{N}_{[1,m]}$. The updated schedule at $k+1$ is
\begin{equation*}
	\tau_j[k+1] = \tau_j[k]+\gamma p_j[k], \quad \forall j \in \mathbb{N}_{[1,m]}.
\end{equation*}
Note that $\tau_j[k+1]\in \mathcal{D}_j$ since \eqref{eqn:local_Quad_approxc} must hold and through the line search condition \eqref{eq:SQPFeasCon} then
\begin{equation*}
	g_i^*((\tau_j[k]+\gamma p_j[k])_{j=1}^m) = g_i^*((\tau_j[k+1])_{j=1}^m) \leq 0.
\end{equation*}
Hence $(\tau_j[k+1])_{j=1}^m$ satisfies \eqref{eqn:Feasble_Q}. As $(\tau_j[0])_{j=1}^m$ is assumed to be feasible, statement one follows by induction.

% \subsubsection{Second statement - termination and cost improvement}
{\it ii)}
Since $B[0]$ is postive definite, $B[k]$ is positive definite for all $k \in \mathbb{N}_0$; see \cite[Chapter 18]{Nocedal2006}. Let
\begin{equation*}
	m_k(p) := \frac{1}{2}p^\top B[k] p + \sum_{j = 1}^m p_j \frac{d}{d \tau_j} h_j(\tau_j[k] ) + \sum_{j = 1}^m h_j(\tau_j[k]).
\end{equation*}
Since $B[k]$ is positive definite, for $p\neq 0$,
\begin{equation}\label{eq:mkPostiveDef}
	m_k(p) > \sum_{j = 1}^m p_j \frac{d}{d \tau_j} h_j(\tau_j[k] ) + \sum_{j = 1}^m h_j(\tau_j[k]).
\end{equation}
The optimal objective value for \eqref{eqn:local_Quad_approx} is given by
\begin{equation}\label{eq:fStarMk}
	\tilde{f}^{*}[k] = m_k(p^*[k]).
\end{equation}
An upper bound for $\tilde{f}^{*}[k]$ is $m_k(0)$, because $p=0$ is a feasible solution for \eqref{eqn:local_Quad_approx}.
% Since $p=0$ is a feasible solution to \eqref{eqn:local_Quad_approx}, $m_k(0)$ is a valid upper bound for $f^{*}[k]$. 
 Combining this with \eqref{eq:fStarMk} and \eqref{eq:mkPostiveDef} gives
\begin{align*}
\sum_{j = 1}^m p^*_j[k] \frac{d}{d \tau_j} h_j(\tau_j[k] ) + \sum_{j = 1}^m h_j(\tau_j[k])  
&< \tilde{f}^{*}[k] 
\\
&\leq \sum_{j = 1}^m h_j(\tau_j[k]).
\end{align*}
Thus,
$ \sum_{j = 1}^m p_j^*[k] \frac{d}{d \tau_j} h_j(\tau_j[k] ) < 0$ for $p^*[k] \neq 0$. 

The last inequality, the line-search condition \eqref{eq:LineSearchCon2}, and the fact that $\gamma$ and $\eta$ are positive, combine to give
\begin{align*}
    \sum_{j=1}^m h_j(\tau_j[k]+\gamma p_j^*[k]) &= \sum_{j=1}^m h_j(\tau_j[k+1])  < \sum_{j=1}^m h_j(\tau_j[k]),% \quad \\
    % &\forall p_j[k] \neq 0
\end{align*}
provided $p^*[k]\neq 0$.
This is equivalent to
\begin{equation}
    m_{k+1}(0)  < m_k(0) \label{eqn:MkIneq}
\end{equation}
when $p^*[k]\neq 0$.
Consider first that $p^*[k]=0$ for some $k\neq 0$. Then the termination condition \eqref{eq:TerminationCond} is satisfied and statement two holds by \eqref{eqn:MkIneq} and the fact $k>0$. Otherwise, $p^*[k]\neq 0$ for all $k \in \mathbb{N}\cup\{0\}$ by hypothesis. In this case, the sequence $(m_k(0))_{k=0}^\infty$ is monotonically decreasing, and it is bounded since $h_j$ are continuously differentiable and each delay must lie in a closed bounded set $\mathcal{D}_j$.  As a result, this sequence converges to a value denoted $\underline{m}$ i.e., $\lim_{k\to \infty} m_k(0) = \underline{m} < m_0(0)$. This implies that for any $\epsilon > 0$ there exists a positive integer $K$ such that 
\begin{equation}\label{eqn:mkCauchy}
	m_{k-1}(0) - m_{k}(0) < \epsilon, \quad  k \geq K.
\end{equation}
Let $\epsilon < -\frac{\epsilon_{step}}{\Gamma[k]} \eta \sum_{j=1}^m p_j^*[k-1] \frac{d}{d \tau_j} h_j(\tau_j[k-1])$. Then combining the line-search condition \eqref{eq:LineSearchCon2} and \eqref{eqn:mkCauchy} leads to
\begin{equation*}
    \gamma < \frac{\epsilon^{s}}{\Gamma[k]}, \quad k \geq K.
\end{equation*}
Consequently, $\gamma \|p^*[k-1]\|_{\infty} < \epsilon^{s}, ~k \geq K,$ 
	as $\|p^*[k-1]\|_{\infty}\leq \Gamma[k]$.
Hence, the SQPA terminates after $K\geq 1$ steps and $m_{K}(0)<m_0(0)$.

\end{document}